\documentclass[aos,MSNbibl,seceqn,dvips]{arximspdf}
\usepackage[linesnumbered,vlined,ruled]{algorithm2e}
\usepackage{dcolumn}
\usepackage{graphicx}

\doi{10.1214/13-AOS1097} 
\volume{41}
\issue{2}
\pubyear{2013}
\firstpage{772}
\lastpage{801}

\makeatletter

\newcolumntype{d}[1]{D{.}{.}{#1}}

\newcommand{\rright}{\right}
\newcommand{\lleft}{\left}

\newtheorem{theorem}{Theorem}[section]
\newtheorem{corollary}{Corollary}[section]
\newtheorem{proposition}{Proposition}[section]
\newtheorem{lemma}{Lemma}[section]

\newtheorem{conditiongr}{Condition}
\newtheorem{conditionsp}{Condition}

\newproclaim{remark}{Remark}[section]
\newproclaim{example}{Example}
\newproclaim{definition}{Definition}

\newcommand{\RR}{\mathbb{R}}
\newcommand{\hf}{{1/2}}
\newcommand{\eps}{\varepsilon}
\newcommand{\Ex}{\mathbb{E}}
\newcommand{\supp}{\operatorname{supp}}

\newcommand{\pn}{p_n}
\newcommand{\class}{\mathcal{F}}
\newcommand{\nspike}{\bar{m}}

\newcommand{\subsp}{\mathcal{P}}
\newcommand{\thr}{\gamma}
\newcommand{\cangle}{{\theta}}
\newcommand{\orac}{\mathrm{o}}

\newcommand{\eval}{\ell}
\newcommand{\Ks}{K_s}

\newcommand{\indi}[1]{\mathrm{1}_{(#1)}}

\makeatother

\begin{document}
\begin{frontmatter}

\title{Sparse principal component analysis and iterative~thresholding}
\runtitle{Sparse principal component analysis}

\begin{aug}
\author[A]{\fnms{Zongming} \snm{Ma}\corref{}\ead[label=e1]{zongming@wharton.upenn.edu}}
\runauthor{Z. Ma}
\affiliation{University of Pennsylvania}
\address[A]{Department of Statistics\\
The Wharton School\\
University of Pennsylvania\\
Philadelphia, Pennsylvania 19104\\
USA\\
\printead{e1}} 
\end{aug}

\received{\smonth{9} \syear{2012}}
\revised{\smonth{1} \syear{2013}}

%
\begin{abstract}
Principal component analysis (PCA) is a classical dimension reduction
method which projects data onto the principal subspace spanned by the
leading eigenvectors of the covariance matrix. However, it behaves
poorly when the number of features $p$ is comparable to, or even much
larger than, the sample size $n$. In this paper, we propose a new
iterative thresholding approach for estimating principal subspaces in
the setting where the leading eigenvectors are sparse.
Under a spiked covariance model, we find that the new approach recovers
the principal subspace and leading eigenvectors consistently, and even
optimally, in a range of high-dimensional sparse settings. Simulated
examples also demonstrate its competitive performance.
\end{abstract}

%

\end{frontmatter}

\section{Introduction}
\label{secintro}

In many contemporary datasets, if we organize the $p$-dimensional observations
$x_1,\ldots, x_n$,
into the rows of an $n\times p$ data matrix $X$, the number of features
$p$ is often comparable to, or even much larger than, the sample size
$n$. For example, in biomedical studies, we usually have measurements
on the expression levels of tens of thousands of genes, but only for
tens or hundreds of individuals. One of the crucial issues in the
analysis of such ``large $p$'' datasets is dimension reduction of the
feature space.

As a classical method, principal component analysis (PCA)
\cite{pear01,hote33}
reduces dimensionality by projecting the data onto the \textit{principal
subspace} spanned by the $m$ leading eigenvectors of the
population covariance matrix $\Sigma$, 
which represent the principal modes of variation.
In principle, one expects that for some $m < p$, most of the variance
in the data is captured by these $m$ modes. Thus, PCA reduces the
dimensionality of the feature space while retaining most of the
information in data.
In addition, projection to a low-dimensional space enables
visualization of the data. In practice, $\Sigma$ is unknown. Classical
PCA then estimates the leading population eigenvectors by those of the
sample covariance matrix $S$. It performs well in the traditional data
setting where $p$ is small and $n$ is large~\cite{ande63}.


In high-dimensional settings, a collection of data can be modeled by a
low-rank signal plus noise structure, and PCA can be used to recover
the low-rank signal. In particular, each observation vector $x_i$ can
be viewed as an independent instantiation of the following generative model:
%
%
\begin{equation}
\label{eqfactor-model} x_i = \mu+ Au_i + \sigma
z_i.
\end{equation}
Here, $\mu$ is the mean vector, $A$ is a $p\times\nspike$
deterministic matrix of factor loadings, $u_i$ is an $\nspike$-vector
of random factors, $\sigma>0$ is the noise level and $z_i$ is a
$p$-vector of white noise.
For instance, in chemometrics, $x_i$ can be a vector of the logarithm
of the absorbance or reflectance spectra measured with noise, where the
columns of $A$ are characteristic spectral responses of different
chemical components, and $u_i$'s the concentration levels of these
components~\cite{vafi09}. 
The number of observations are relatively few compared with the number
of frequencies at which the spectra are measured. In econometrics,
$x_i$ can be the returns for a collection of assets, where the $u_i$'s
are the unobservable random factors~\cite{tsay05}. 
The assumption of additive white noise is reasonable for asset returns
with low frequencies (e.g., monthly returns of stocks).
Here, people usually look at tens or hundreds of assets simultaneously,
while the number of observations are also at the scale of tens or
hundreds. In addition, model (\ref{eqfactor-model}) represents a big
class of signal processing problems~\cite{wk85}. 
Without loss of generality, we assume $\mu= 0$ from now on.

In this paper, our primary interest lies in PCA of high-dimensional
data generated as in (\ref{eqfactor-model}).
Let the covariance matrix of $u_i$ be $\Phi$ which is of full rank.
Suppose that $A$ has full column rank and that $u_i$ and $z_i$ are independent.
Then the covariance matrix of $x_i$ becomes
%
%
\begin{equation}
\label{eqSigma-spike} \Sigma= A\Phi A' + \sigma^2 I =
\sum_{j=1}^{\nspike} \lambda_j^2
q_j q_j' + \sigma^2 I.
\end{equation}
Here, $\lambda_1^2\geq\cdots\geq\lambda_{\nspike}^2 > 0$ are the
eigenvalues of $A\Phi A'$, with $q_j$, $j=1,\ldots, \nspike$, the
associated eigenvectors. Therefore, the $j$th eigenvalue of $\Sigma$
is $\lambda_j^2+\sigma^2$ for $j=1,\ldots, \nspike$, and $\sigma^2$
otherwise. Since there are $\nspike$ spikes $(\lambda_1^2,\ldots,
\lambda_{\nspike}^2)$ in the spectrum of $\Sigma$, (\ref
{eqSigma-spike}) has been called the \textit{spiked covariance model} in
the literature~\cite{john01}.
Note that we use $\lambda_j^2$ to denote the spikes rather than
$\lambda_j$ used previously in the literature~\cite{pajo07}.
For data with such a covariance structure, it makes sense to project
the data onto the low-dimensional subspaces spanned by the first few $q_j$'s.
Here and after, $\nspike$ denotes the number of spikes in the model,
and $m$ is the target dimension of the principal subspace to be
estimated, which is no greater than $\nspike$.

Classical PCA encounters both practical and theoretical difficulties in
high dimensions.\vadjust{\goodbreak}
On the practical side, the eigenvectors found by classical PCA involve
all the $p$ features, which makes their interpretation challenging.
On the theoretical side, the sample eigenvectors are no longer always
consistent estimators. 
Sometimes, they can even be nearly orthogonal to the target direction.
When both $n,p\to\infty$ with $n/p\to c\in(0,\infty)$, at different
levels of rigor and generality, this phenomenon has been examined by a
number of authors~\cite{reva96,lu02,hora04,paul07,nadl08,onat09}
under model (\ref{eqSigma-spike}). See~\cite{juma09} for similar results
when $p\to\infty$ and $n$ is fixed.

In recent years, to facilitate interpretation, researchers have started
to develop sparse PCA methodologies, where they seek a set of sparse
vectors spanning the low-dimensional subspace that explains most of the
variance. See, for example, \mbox{\cite
{jotrud03,zohati06,dagh07,shhu08,solo08,witiha09}}. These approaches
typically start with a certain
optimization formulation of PCA and then induce a sparse solution by
introducing appropriate penalties or constraints.

On the other hand, when $\Sigma$ indeed has sparse leading eigenvectors
in the current basis
(perhaps after transforming the data),
it becomes possible to estimate them consistently under
high-dimensional settings via new estimation schemes.
For example, under normality assumption, when $\Sigma$ only has a
single spike, that is, when $\nspike= 1$ in (\ref{eqSigma-spike}),
Johnstone and Lu~\cite{jolu09} proved consistency of PCA obtained on a
subset of features with large sample variances when the leading
eigenvalue is fixed and $(\log{p})/n\to0$. Under the same single
spike model, if in addition the leading eigenvector has exactly $k$
nonzero loadings, Amini and Wainwright~\cite{amwa09} studied
conditions for recovering the nonzero locations using the methods in
\cite{jolu09} and~\cite{dagh07}, and Shen et al.~\cite{shshma11}
established conditions for consistency of a sparse PCA method in \cite
{shhu08} when $p\to\infty$ and $n$ is fixed. For the more general
multiple component case, Paul and Johnstone~\cite{pajo07} proposed an
augmented sparse PCA method for estimating each of the leading
eigenvectors, and showed that their procedure attains near optimal rate
of convergence under a range of high-dimensional sparse settings when
the leading eigenvalues are comparable and well separated. Notably,
these methods all focus on estimating individual eigenvectors.

In this paper, we focus primarily on finding \textit{principal
subspaces} of $\Sigma$ spanned by sparse leading eigenvectors, as
opposed to finding each sparse vector individually.
One of the reasons is that individual eigenvectors are not identifiable
when some leading eigenvalues are identical or close to each other.
Moreover, if we view PCA as a dimension reduction technique, it is the
low-dimensional subspace onto which we project data that is of the
greatest interest.

We propose a new iterative thresholding algorithm to estimate principal
subspaces, which is motivated by the orthogonal iteration method in
matrix computation. 
In addition to the usual orthogonal iteration steps, an additional
thresholding step is added to seek 
sparse basis vectors
for the subspace.
When $\Sigma$ follows the spiked covariance model and the sparsity of
the leading eigenvectors are characterized by the weak-$\ell_r$
condition (\ref{eqweak-lr}),
the algorithm leads to 
a consistent subspace estimator adaptively over a wide range of
high-dimensional sparse settings, and the rates of convergence are
derived under an appropriate loss function (\ref{eqloss}).
Moreover, for any individual leading eigenvector whose eigenvalue is
well separated from the rest of the spectrum, our algorithm also yields
an eigenvector estimator which adaptively attains optimal rate of
convergence derived in~\cite{pajo07} up to a multiplicative log factor.
In addition, it has appealing model selection property in the sense
that the resulting estimator only involves coordinates with large
signal-to-noise ratios.

The contribution of the current paper is threefold.
First, we propose to estimate principal subspaces. This is natural for
the purpose of dimension reduction and visualization, and avoids the
identifiability issue for individual eigenvectors.
Second, we construct a new algorithm to estimate the subspaces, which
is efficient in computation and easy to implement.
Last but not least,
we derive convergence rates of the resulting estimator under the spiked
covariance model
when the eigenvectors 
are sparse.

The rest of the paper is organized as follows.
In Section~\ref{secalg}, we frame the principal subspace estimation
problem and propose the iterative thresholding algorithm. The
statistical properties and computational complexity of the algorithm
are examined in Sections~\ref{secopt} and~\ref{seccomplex} under
normality assumption. Simulation results in Section~\ref{secnumeric}
demonstrate its competitive performance.
Section~\ref{secproof} presents the proof of the main theorems.

\textit{Reproducible code}: The \textsc{Matlab} package \texttt
{SPCALab} implementing the proposed method and producing the tables and
figures of the current paper is available at the author's website.

\section{Methodology}
\label{secalg}


\subsection{Notation}
We say $x$ is a $p$-vector if $x\in\RR^p$,
and we use $\|x\|_{2}$ to denote its Euclidean norm.
For an $m\times n$ matrix $A$,
its submatrix with rows indexed by $I$ and columns indexed by $J$ is
denoted by $A_{I J}$. If $I$ or $J$ includes all the indices, we
replace it with a dot. For example, $A_{I \cdot}$ is the
submatrix of $A$ with rows in $I$ and all columns.
The spectral norm of $A$ is
$\|A\| = \max_{\|x\|_{2}=1}\|Ax\|_{2}$, 
and the range, that is, the column subspace, of $A$ is $\operatorname{ran}(A)$.
If $m\geq n$, and the columns of $A$ form an orthonormal set in $\RR
^m$, we say $A$ is orthonormal. 

We use $C, C_0, C_1$, etc. to represent constants, though their values
might differ at different occurrences. For real numbers $a$ and $b$,
let $a\vee b = \max(a,b)$ and $a\wedge b = \min(a,b)$.
We write $a_n = O(b_n)$, if there is a constant $C$, such that
$|a_n|\leq Cb_n$ for all $n$, and $a_n = o(b_n)$ if $a_n/b_n \to0$ as
$n\to\infty$.
Moreover, we write $a_n\asymp b_n$ if $a_n = O(b_n)$ and $b_n = O(a_n)$.
Throughout the paper, we use $\nu$ as the generic index for features,
$i$ for observations, $j$ for eigenvalues and eigenvectors and $k$ for
iterations in the algorithm to be proposed.

\subsection{Framing the problem: Principal subspace estimation}

When the covariance matrix $\Sigma$ follows 
model (\ref{eqSigma-spike}), its $j$th largest eigenvalue $\eval
_{j}(\Sigma) = \lambda_j^2 + \sigma^2$ for $j=1,\ldots, \nspike$
and equals $\sigma^2$ for all $j > \nspike$. Let $\operatorname
{span}\{\cdot\}$
denote the linear subspace spanned by the vectors in the curly brackets.
If for some $m\leq\nspike$, $\eval_{m}(\Sigma) > \eval
_{m+1}(\Sigma)$,
the \textit{principal subspace}
\[
\subsp_m = \operatorname{span}\{q_1,\ldots,q_m \}
\]
is defined, regardless of the behavior of the other $\eval
_{j}(\Sigma)$'s.
Therefore, it is an identifiable object for the purpose of estimation.
Note that $\subsp_{\nspike}$ is always identifiable, because $\eval
_{\nspike}(\Sigma) > \eval_{\nspike+1}(\Sigma)$.
The primary goal of this paper is to estimate the principal subspace
$\subsp_m$, for some $m\leq\nspike$
with $\ell_m(\Sigma) > \ell_{m+1}(\Sigma)$.
More precisely, we require
the gap $\ell_m(\Sigma) - \ell_{m+1}(\Sigma) = \lambda_m^2-
\lambda_{m+1}^2\asymp\lambda_1^2$. Note that such an $m$ always
exists, for example, the largest $m\leq\nspike$ such that $\lambda
_m^2\asymp\lambda_1^2$.
We allow the case of $m<\nspike$ partly because under certain
circumstances, one might not be interested in $\subsp_{\nspike}$
directly. For example, to visualize the data, one might want to estimate
$\subsp_2$ or $\subsp_3$ while $\nspike$ could be larger than $3$.
In addition, sometimes $\subsp_{\nspike}$ might not be consistently
estimable while some smaller principal subspace $\subsp_m$ is.
In most part of the paper, we assume that an appropriate $m$ is given
for convenience.
In Section~\ref{secrank}, we discuss how to choose $m$ and how to
estimate $\bar{m}$ under normality assumption.


To measure the accuracy of an estimator $\widehat{\mathcal{S}}$ for a
subspace $\mathcal{S}$,
note that each linear subspace is associated with a unique
projection matrix onto it.
Let $P$ and $\widehat{P}$ be the projection matrices associated with
$\mathcal{S}$ and $\widehat{\mathcal{S}}$, respectively. The
distance between $\mathcal{S}$ and $\widehat{\mathcal{S}}$\vspace*{2pt}
is given
by the spectral norm of the difference between $P$ and $\widehat{P}$:
$\operatorname{dist}(\mathcal{S}, \widehat{\mathcal{S}}) = \|P -
\widehat
{P}\|$; see 
\cite{gova96}, Section 2.6.3.
Thus, we can define a loss function by the squared distances between
$\mathcal{S}$ and~$\widehat{\mathcal{S}}$,
%
%
\begin{equation}
\label{eqloss} L(\mathcal{S}, \widehat{\mathcal{S}}) = \operatorname{dist}^2(
\mathcal{S}, \widehat{\mathcal{S}}) = \|P - \widehat{P}\|^2.
\end{equation}
By definition, this loss function measures the maximum possible
discrepancy between the projections of any unit vector onto the two
subspaces. The loss ranges in $[0,1]$, and equals zero if and only if
$\widehat{\mathcal{S}} = \mathcal{S}$.
When $\operatorname{dim}(\widehat{\mathcal{S}}) \neq
\operatorname{dim}(\mathcal{S})$, we have $L(\mathcal{S},\widehat
{\mathcal{S}})=1$.
Geometrically, it equals the squared sine of the largest canonical
angle between $\mathcal{S}$ and $\widehat{\mathcal{S}}$
(\cite{stsu90}, Theorem 5.5).
Throughout the paper, we use the loss function (\ref{eqloss}) for
principal subspace estimation.


\subsection{Orthogonal iteration} 
\label{secor-it}
Given a positive definite matrix $A$, a standard technique to compute
its leading 
eigenspace is orthogonal iteration~\cite{gova96}.
When only the first eigenvector is sought, it is also known as the
power method.

To state the orthogonal iteration method, we note that for any
$p\times m$ matrix $T$, when $p\geq m$, we could decompose it into the
product of two matrices $T = QR$, where $Q$ is $p\times m$ orthonormal
and $R$ is $m\times m$ upper triangular. This decomposition is called\vadjust{\goodbreak}
QR factorization and can be computed using Gram--Schmidt
orthogonalization and other numerical methods~\cite{gova96}.
Suppose $A$ is $p\times p$, and we want to compute its leading
eigenspace of dimension $m$.
Starting with a $p\times m$ orthonormal matrix $Q^{(0)}$,
orthogonal iteration generates a sequence of $p\times m$ orthonormal
matrices $Q^{(k)}$, $k=1,2,\ldots\,$, by alternating the following
two steps till convergence:
\begin{longlist}[(2)]
\item[(1)] Multiplication: $T^{(k)} = AQ^{(k-1)}$;
\item[(2)] QR factorization: $Q^{(k)}R^{(k)} = T^{(k)}$.
\end{longlist}
Denote the orthonormal matrix at convergence by $Q^{(\infty)}$.
Then its columns are the leading eigenvectors
of $A$, and $\operatorname{ran}(Q^{(\infty)})$ gives the eigenspace. In
practice, one terminates the iteration once $\operatorname
{ran}(Q^{(k)})$ stabilizes.

When we apply orthogonal iteration directly to the sample covariance
matrix $S$, it gives the classical PCA result, which could be
problematic in high dimensions. Observe that all the $p$ features are
included in orthogonal iteration. When the dimensionality is high,
not only the interpretation is hard, but the variance accumulated
across all the features becomes so high that it makes consistent
estimation impossible.

If the eigenvectors spanning $\subsp_m$ are sparse in the current
basis, one sensible way to reduce estimation error is to focus only on
those features at which the leading eigenvectors have large values,
and to estimate other features by zeros. Of course, one introduces
bias this way, but hopefully it is much smaller compared to the amount
of variance thus reduced.

The above heuristics
lead to the estimation scheme in the next subsection which
incorporates this feature screening idea in orthogonal iteration.

\subsection{Iterative thresholding algorithm}
\label{secit-th}

Let $S = \frac{1}{n}\sum_{i=1}^n x_i x_i'$ be the sample covariance matrix.
An effective way to incorporate feature screening into orthogonal
iteration is to ``kill'' small coordinates of the $T^{(k)}$
matrix after each multiplication step, which leads to the estimation
scheme summarized in Algorithm~\ref{algspcait}. Although the later
theoretical study is conducted under normality assumption, Algorithm
\ref{algspcait} itself is not confined to normal data.


\begin{algorithm}\label{algspcait}
\caption{ITSPCA (Iterative thresholding sparse PCA)}
{\fontsize{10pt}{12pt}\selectfont{\KwIn{
\begin{enumerate}[(4)]
\item[\hphantom{\textbf{In}}(1)] Sample covariance matrix ${S}$;
\item[\hphantom{\textbf{In}}(2)] Target subspace dimension $m$;
\item[\hphantom{\textbf{In}}(3)] Thresholding function $\eta$, and threshold levels $\thr
_{nj}$, $j=1,\ldots, m$;
\item[\hphantom{\textbf{In}}(4)] Initial orthonormal matrix $\widehat{Q}{}^{(0)}$.
\end{enumerate}
}
\KwOut{Subspace estimator $\widehat{\subsp}_m =
\operatorname{ran}(\widehat{Q}{}^{(\infty)})$, where $\widehat
{Q}{}^{(\infty)}$ denotes the $\widehat{Q}{}^{(k)}$ matrix at
convergence. }
\Repeat{convergence}{
Multiplication: $T^{(k)} = (t^{(k)}_{\nu j}) =
{S}\widehat{Q}{}^{(k-1)}$;

Thresholding: $\widehat{T}{}^{(k)} = (\widehat{t}{}^{(k)}_{\nu j})$, with
$\widehat{t}{}^{(k)}_{\nu j}= \eta(t^{(k)}_{\nu j}, \thr_{nj})$;

QR factorization: $\widehat{Q}{}^{(k)}\widehat{R}{}^{(k)} =
\widehat{T}{}^{(k)}$;
}}}
\end{algorithm}
%


In addition to the two basic orthogonal iteration steps, Algorithm \ref
{algspcait} adds a thresholding step in between them, where we
threshold each element of $T^{(k)}$ with a user-specified
thresholding function $\eta$ which satisfies
%
%
\begin{equation}
\label{eqeta}\qquad
\bigl|\eta(t,\thr) - t\bigr|\leq\thr\quad\mbox{and}\quad\eta
(t,\thr)\indi
{|t|\leq
\thr} = 0\qquad\mbox{for all $t$ and all $\thr> 0$}.
\end{equation}
Here, $\indi{E}$ denotes the indicator function of an event $E$. We
note that both hard-thresholding
$\eta_H(t,\thr) = t \indi{|t|>\thr}$
and soft-thresholding
$\eta_S(t,\thr) = \operatorname{sgn}(t)(|t|-\thr)_+$
satisfy (\ref{eqeta}).
So does any $\eta$ sandwiched by them,
such as that resulting from a SCAD criterion~\cite{fali01}.
In $\eta(t,\thr)$, the parameter $\thr$ is called the threshold
level. In Algorithm~\ref{algspcait}, for each column of $T^{(k)}$, a
common threshold level $\thr_{nj}$ needs to be specified for
all its elements, which remains unchanged across iterations. The
subscripts of $\thr_{nj}$ indicate that it depends on both the size of
the problem $n$ and the index $j$ of the column it is applied to.

%
%
\begin{remark}
The ranges of $\widehat{Q}{}^{(k)}$ and $\widehat{T}{}^{(k)}$
are the same because QR factorization only amounts to a basis change
within the same subspace.
However, as in orthogonal iteration, the QR step is essential for
numerical stability, and should not be omitted.
Moreover, although the algorithm is designed for subspace estimation,
the column vectors of $\widehat{Q}{}^{(\infty)}$ can be used as
estimators of leading eigenvectors.
\end{remark}

\subsubsection*{Initialization}
Algorithm~\ref{algspcait} requires an initial orthonormal matrix
$\widehat{Q}{}^{(0)}$.
It can be generated from the ``diagonal thresholding'' sparse PCA algorithm
\cite{jolu09}.
Its multiple eigenvector version is summarized in Algorithm~\ref{algdtspca}.
Here,\vspace*{1pt} for any set $I$, $\operatorname{card}(I)$ denotes its cardinality.
\begin{algorithm}\label{algdtspca}
\caption{DTSPCA (Diagonal thresholding sparse PCA)}
{\fontsize{10pt}{12pt}\selectfont{
\KwIn{
\begin{enumerate}
\setlength{\itemsep}{0.15em}
\setlength{\parskip}{0em}
{
\item[(1)] Sample covariance matrix ${S}$;
%
\item[(2)] Diagonal thresholding parameter $\alpha_n$.
}
\end{enumerate}
}

\KwOut{Orthonormal matrix {$\widehat{Q}_B$}.}

Variance selection: select the set $B$ of coordinates (which are
likely to have ``big'' signals):
\[
B = \bigl\{\nu\dvtx s_{\nu\nu} \geq\sigma^2 (1 +
\alpha_n)\bigr\}; 
\]

Reduced PCA: compute the eigenvectors, $\widehat{q}_1^B,\ldots,
\widehat{q}_{\operatorname{card}(B)}^B$, of the submatrix $S_{BB}$;

Zero-padding: construct $\widehat{Q}_B =
[\widehat{q}_1,\ldots,\widehat{q}_{\operatorname{card}(B)}]$ such that
\begin{eqnarray*}
\widehat{q}_{jB} &=& \widehat{q}_j^B,\qquad \widehat
{q}_{jB^c} = 0,\qquad j=1,\ldots, \operatorname{card}(B).\\[-28pt]
\end{eqnarray*}
}}
\end{algorithm}
Given the output $\widehat{Q}_B = [\widehat{q}_1,\ldots,\widehat
{q}_{\operatorname{card}(B)}]$ of Algorithm~\ref{algdtspca}, we take
$\widehat
{Q}{}^{(0)} = [\widehat{q}_1,\ldots,\widehat{q}_m]$.
When\vspace*{1pt} $\sigma^2$ is unknown, we could replace it by an estimator
$\widehat\sigma^2$ in the definition of $B$.
For example, for normal data, Johnstone and Lu~\cite{jolu09} suggested
%
%
\begin{equation}
\label{eqhat-sigma} \widehat\sigma^2 = \operatorname{median}\Biggl(
\frac{1}{n} {\sum_{i=1}^n}
x_{i\nu}^2 \Biggr).
\end{equation}
%
When available, subject knowledge could also be incorporated into the
construction of $\widehat{Q}{}^{(0)}$.
Algorithm~\ref{algspcait} also requires inputs for the $\thr_{nj}$'s
and subspace dimension~$m$. Under normality assumption, we give
explicit specification for them in (\ref{eqthr-nj}) and (\ref{eqm}) later.
Under the conditions of the later Section~\ref{secopt},
$B$ is nonempty with probability tending to $1$, and so $\widehat
{Q}{}^{(0)}$ is well defined.

\subsubsection*{Convergence}
For normal data,
to obtain the error rates in later Theorems~\ref{thmspecial} and \ref
{thmerror},
we can terminate Algorithm~\ref{algspcait} after $\Ks$ iterations
with $\Ks$ given in (\ref{eqstop}).\vadjust{\goodbreak}
In practice, one could also stop iterating if the difference between
successive iterates becomes sufficiently small, for example, when
$L(\operatorname{ran}(\widehat{Q}{}^{(k)}),\operatorname{ran}(\widehat
{Q}{}^{(k+1)}))
\leq n^{-2}$.
We suggest this empirical stopping rule because $n^{-2}$
typically tends to zero faster than the rates we shall obtain,
and so intuitively it should not change the statistical performance of
the resulting estimator.
In simulation studies reported in Section~\ref{secnumeric}, the
difference in numerical performance between the outputs based on this
empirical stopping rule and those based on the theoretical rule (\ref
{eqstop}) is negligible compared to the estimation errors.
Whether Algorithm~\ref{algspcait} always converges numerically
is an interesting question left for possible future research.

\subsubsection*{Bibliographical note}
When $m=1$, Algorithm~\ref{algspcait} is similar to the algorithms
proposed in 
\cite{shhu08,witiha09}
and~\cite{yuzh11}.
When $m > 1$, all these methods
propose to iteratively find the first leading eigenvectors of residual
covariance matrices, which becomes different from our approach.

\section{Statistical properties}
\label{secopt}
This section is devoted to analyzing the statistical properties of
Algorithm~\ref{algspcait} under normality assumption.
After some preliminaries, we first establish the convergence rates for
subspace estimation in a special yet interesting case in Section \ref
{secspecial}. Then we introduce a set of general assumptions in
Section~\ref{secmodel} and a few key quantities in Section \ref
{secprelim}. Section~\ref{secresult} states the main results, which include
convergence rates for principal subspace estimation under general
assumptions and a correct exclusion property.
In addition, we derive rates for estimating individual eigenvectors.
For conciseness, we first state all the results assuming a suitable
target subspace dimension $m\leq\nspike$ is given. In Section \ref
{secrank}, we discuss how to choose $m$ and estimate $\nspike$ based
on data.

We start with some preliminaries. Under normality assumption,
$x_1,\ldots,x_n$ are i.i.d. $N_p(0,\Sigma)$ distributed, with $\Sigma
$ following 
model (\ref{eqSigma-spike}).
Further assume $\sigma^2$ is known---though this assumption could be
removed by estimating $\sigma^2$ using, say, $\widehat\sigma^2$ in
(\ref{eqhat-sigma}).
Since one can always scale the data first, we assume $\sigma^2 = 1$
from now on.
Thus, (\ref{eqfactor-model}) reduces to the 
orthogonal factor
form
%
%
\begin{equation}
\label{eqmodel} x_i = \sum_{j=1}^{\nspike}
\lambda_j v_{ij} q_j + z_i,\qquad i=1,\ldots, n.
\end{equation}
Here, 
$v_{ij}$ are i.i.d. standard normal random factors, which are
independent of the i.i.d. white noise vectors $z_i\sim N_p(0,I)$, and
$\{q_j, 1\leq j\leq\nspike\}$ is a set of leading eigenvectors of
$\Sigma$.
In what follows, we use $n$ to index the size of the problem. So the
dimension $p = p(n)$ and the spikes $\lambda_j^2 = \lambda_j^2(n)$
can be regarded as functions of $n$, while both $\nspike$ and $m$
remain fixed as $n$ grows.

Let $\pn= p\vee n$.
We obtain the initial matrix $\widehat{Q}{}^{(0)}$ in
Algorithm~\ref{algspcait} by applying Algorithm~\ref{algdtspca} with
%
%
\begin{equation}
\label{eqalpha-n} \alpha_n = \alpha\biggl[\frac{\log(\pn)}{n}
\biggr]^{\hf}.
\end{equation}
%
In Algorithm~\ref{algspcait}, the threshold levels are set at
%
%
\begin{equation}
\label{eqthr-nj} 
\thr_{nj} = \thr\biggl[
\ell^B_{j} \frac{\log(\pn)}{n} \biggr]^{\hf},\qquad j=1,\ldots,m.
\end{equation}
Here, $\alpha$ and $\thr$ are user specified constants, and {$\ell^B
_j = \eval_j(S_{BB})\vee1$} with $\eval_j(S_{BB})$ the $j$th largest
eigenvalue of $S_{BB}$, where the set $B$ is obtained in step 1 of
Algorithm~\ref{algdtspca}.
For theoretical study, we always stop Algorithm~\ref{algspcait} after
$\Ks$ iterations, where for $h(x) = x^2/(x+1)$,
%
%
\begin{equation}
\label{eqstop} \Ks= \frac{1.1 \cdot\ell^B_1}{\ell^B_{m}- \ell^B_{m+1}}
\biggl[ \biggl(1 + \frac{1}{\log{2}}
\biggr)\log{n} + 0\vee\log h\bigl(\ell^B _1-1\bigr)
\biggr].
\end{equation}
Under the conditions of Theorem~\ref{thmspecial} or of Theorems \ref
{thmerror} and~\ref{thmexclude}, or when $m$ is defined by (\ref
{eqm}), we have $\ell^B_m\neq\ell^B_{m+1}$ and $\Ks< \infty$ with
probability $1$.

\subsection{A special case}
\label{secspecial}
To facilitate understanding, we first state the convergence rates for
principal subspace estimation in a special case.

Consider the asymptotic setting where $n\to\infty$ with $p\geq n$ and
$(\log{p}) / n\to0$, while 
the spikes $\lambda_1^2 \geq\cdots\geq\lambda_m^2 > \lambda
_{m+1}^2 \geq\cdots\geq\lambda_{\nspike}^2 > 0$ remain unchanged.
Suppose that the $q_j$'s are sparse in the sense that, for some $r\in
(0,2)$, the $\ell_r$ norm of the eigenvectors are uniformly bounded by
$s$, that is,
$\|q_j\|_{r} = (\sum_{\nu=1}^p |q_{\nu j}|^r )^{1/r} \leq s$,
for $j=1,\ldots,\nspike$,
where $s\geq1$ is an absolute constant.

Recall that $h(x) = x^2/(x+1)$.
Under the above setup, we have the following upper bound
for subspace estimation error.
%
%
\begin{theorem}
\label{thmspecial}
Under the above setup,
for sufficiently large constants $\alpha, \thr{> 2\sqrt{3}}$ in
(\ref{eqalpha-n}) and (\ref{eqthr-nj}), there exist constants
$C_0$, $C_1 = C_1(\thr, r, m)$ and $C_2$, such that for sufficiently
large $n$,
uniformly over all $\Sigma$ with $\|q_j\|_r\leq s$ for $1\leq j\leq
\nspike$,
with probability at least $1-C_0\pn^{-2}$,
we have $\Ks\asymp\log{n}$ and
the subspace estimator $\widehat{\subsp}_m^{({\Ks})} = \operatorname
{ran}(\widehat{Q}{}^{({\Ks})})$ 
of Algorithm~\ref{algspcait} satisfies
\[
L\bigl(\subsp_m, \widehat{\subsp}_m^{({\Ks})}
\bigr) \leq C_1 \nspike s^r \biggl[ \frac{\log{p}}{ n h(\lambda_m^2) }
\biggr]^{1-r/2} + C_2 g_m(\lambda)
\frac{\log{p}}{n}, 
\]
%
where $g_m(\lambda) = \frac{(\lambda_1^2+1)(\lambda
_{m+1}^2+1)}{(\lambda_m^2-\lambda_{m+1}^2)^2}$.
\end{theorem}

The upper bound in Theorem~\ref{thmspecial} consists of two terms.
The first is a ``nonparametric'' term, which can be decomposed as the
product of two components.
The first component,
$\nspike s^r [nh(\lambda_m^2)/\log{p}]^{r/2}$, up to a multiplicative
constant, bounds the number of coordinates used in estimating the
subspace, while the second component, $\log{p}/[nh(\lambda_m^2)]$,
gives the average error per coordinate.
The second term in the upper bound, $g_m(\lambda)(\log{p})/n$, up to
a logarithmic factor, has the same form as the cross-variance term in
the ``fixed $p$, large $n$'' asymptotic limit for classical PCA; cf.
\cite{ande63}, Theorem 1.
We call it a ``parametric'' error term, because it always arises when
we try to separate the first $m$ eigenvectors from the rest, regardless
of how sparse they are. Under the current setup, both terms 
converge to $0$ as $n\to\infty$, which establishes the consistency of
our estimator.

To better understand the upper bound, we compare it with an existing
lower bound.
Suppose $\lambda_1^2 > \lambda_2^2$. Consider the simplest case where
$m = 1$. Then, estimating $\subsp_1$ is the same as estimating the
first eigenvector $q_1$. For estimating an individual eigenvector
$q_j$, Paul and Johnstone
\cite{pajo07} considered the loss function
$l(q_j, \tilde{q}_j) = \|q_j - \operatorname{sgn}(q_j'\tilde
{q}_j)\tilde{q}_j\|_{2}^2$.
Here, the $\lambda_j^2$'s, $s$ and $r$ are fixed and $p\geq n$, so
when $n$ is large, $s^r [nh(\lambda_1^2)/\log{p}]^{r/2} \leq
Cp^{1-c}$ for some $c\in(0,1)$. For this case, Theorem 2 in \cite
{pajo07} asserts that for any estimator $\widehat{q}_1$,
\[
\sup_{\|q_j\|_r \leq s,\forall j} \Ex l(q_1,
\widehat{q}_1) \geq C_1 s^r \biggl[
\frac{\log{p}}{ n h(\lambda_1^2) } \biggr]^{1-r/2} + C_2 \frac
{g_1(\lambda)}{n}.
\]
Let\vspace*{1pt} $\widehat{\subsp}_1 = \operatorname{span}\{\widehat{q}_1\}$.
We have $\frac{1}{2} l(q_1, \widehat{q}_1) \leq L(\subsp_1, \widehat
\subsp_1) \leq l(q_1, \widehat{q}_1)$. So the above lower bound also
holds for any $\widehat{\subsp}_1$ and $\Ex L(\subsp_1, \widehat
{\subsp}_1)$.
Note that in both Theorem~\ref{thmspecial} and the last display, the
nonparametric term is dominant, and so both the lower and upper bounds
are of order $[(\log p)/n]^{1-r/2}$.
Therefore, Theorem~\ref{thmspecial} shows that the estimator from
Algorithm~\ref{algspcait} is rate optimal.


Since $\alpha_n$ and $\thr_{nj}$ and the stopping rule (\ref
{eqstop}) do not involve any unknown parameter, the theorem
establishes the adaptivity of our estimator: the 
optimal rate of convergence in Theorem~\ref{thmspecial} is obtained
without any knowledge of the power $r$, the radius $s$ or the spikes
$\lambda_j^2$.
Last but not least, 
the estimator could be obtained in $O(\log{n})$ iterations and
holds for all thresholding function $\eta$ satisfying (\ref{eqeta}).


Later in Section~\ref{secresult}, Theorem~\ref{thmerror}
establishes analogous convergence rates, but for a much wider range of
high-dimensional sparse settings.
In particular, the above result will be extended simultaneously along
two directions:
\begin{longlist}[(2)]
\item[(1)] the spikes $\lambda_1^2,\ldots, \lambda_{\nspike}^2$
will be allowed to scale as $n\to\infty$, and
$\lambda_{m+1}^2,\ldots,\break
\lambda_{\nspike}^2$ could even be of smaller order as compared to
the first $m$ spikes;
\item[(2)] each individual eigenvector $q_j$ will be constrained to a
weak-$\ell_r$ ball of radius $s_j$ (which contains the $\ell_r$ ball
of the same radius), and the radii $s_j$'s will be allowed to diverge
as $n\to\infty$.
\end{longlist}


\subsection{Assumptions}
\label{secmodel}
We now state assumptions for the general theoretical results in Section
\ref{secresult}.

As outlined above, the first extension of the special case is to allow
the spikes $\lambda_j^2 = \lambda_j^2(n) > 0$ to change with $n$,
though the dependence will usually not be shown explicitly.
Recall that $p_n = p\vee n$; we impose the following growth rate
condition on $p$ and the $\lambda_j^2$'s.
\renewcommand{\theconditiongr}{GR}
\begin{conditiongr}\label{condGR}
As $n\to\infty$, we have:
\begin{longlist}[(3)]
\item[(1)] the dimension $p$ satisfies $(\log{p})/n = o(1)$;
\item[(2)] the\vspace*{1pt} largest spike $\lambda_1^2$ satisfies
$\lambda_1^2 = O(\pn)$;
the smallest spike $\lambda_{\nspike}^2$ satisfies $\log(\pn) =
o(n\lambda_{\nspike}^4)$; and their ratio satisfies
${\lambda_1^2}/{\lambda_{\nspike}^2} = O(n[\log(\pn)/n]^{\hf+r/4})$;
\item[(3)] $\lim_{n\to\infty} \lambda_1^2/(\lambda_j^2-\lambda
_{j+1}^2)\in[1,\infty]$ exists for $j=1,\ldots, \nspike$, with
$\lambda^2_{\nspike+1}=0$.
\end{longlist}
\end{conditiongr}

The first part of Condition~\ref{condGR} requires the dimension to grow at a
sub-exponential rate of the sample size. The second part ensures that
the spikes 
grow at most at linear rate with $\pn$, 
and are all of larger magnitude than $\sqrt{\log(\pn)/n}$.
In addition, the condition on the ratio $\lambda_1^2/\lambda_{\nspike
}^2$ allows us to deal with the interesting cases where the first
several spikes scale at a faster rate with $n$ than the others. This is
more flexible than the assumption previously made in~\cite{pajo07}
that all the spikes grow at the same rate.
The third part requires $\lim_{n\to\infty} \lambda_1^2/(\lambda
_j^2-\lambda_{j+1}^2)$ to exist for each $1\leq j\leq\nspike$, but
the limit can be infinity.

Turn to the sparsity assumption on the $q_j$'s. We first make a mild
extension from $\ell_r$ ball to weak-$\ell_r$ ball~\cite{dono93}. To
this end, for any $p$-vector $u$, order its coordinates by magnitude as
$|u|_{(1)}\geq\cdots\geq|u|_{(p)}$. We say that $u$ belongs to the
weak-$\ell_r$ ball of radius $s$, denoted by $u\in w\ell_{r}(s)$, if
%
%
\begin{equation}
\label{eqweak-lr} |u|_{(\nu)} \leq
{s} {\nu^{-1/r}}\qquad
\mbox{for all $\nu$}.
\end{equation}
For $r\in(0,2)$, the above condition implies rapid decay of the
ordered coefficients of $u$, and thus describes its sparsity. For
instance, consider $u = (1/\sqrt{k},\ldots, 1/\sqrt{k},0,\ldots,0)'$
with exactly $k$ nonzero entries all equal to $1/\sqrt{k}$. Then, for
fixed $r\in(0,2)$, we have $u\in w\ell_{r}(k^{1/r-\hf})$. In
particular, when $k=1$, $u\in w\ell_r(1)$. Note that weak-$\ell_r$
ball extends $\ell_r$ ball, because $\|u\|_r \leq s$, that is, $u\in
\ell_r(s)$, implies $u\in w\ell_r(s)$.

In what follows, we assume that for some fixed $r\in(0,2)$ and all
$j\leq\nspike$, $q_j\in w\ell_r(s_j)$ for some $s_j > 1$. We choose
to use the notion of ``weak-$\ell_r$ decay,'' because it provides a
unified framework for several different notions of sparsity, which is
convenient for analyzing a statistical estimation problem from a
minimax point of view~\cite{dono93}.
Hence, at any fixed $n$, we will consider whether Algorithm \ref
{algspcait} performs uniformly well on $n$ i.i.d. observations $x_i$
generated by (\ref{eqmodel}) whose covariance matrix $\Sigma$
belongs to the following uniformity class:
\[
\class_n 
= \Biggl\{\Sigma_{p\times p} =
\sum_{j=1}^{\nspike}\lambda_j^2
q_j q_j'+I\dvtx q_j\in w
\ell_r(s_j), \forall j \Biggr\}.
\]

For general results, we allow the radii $s_j$'s to depend on or even
diverge with $n$, though we require that they do not grow too rapidly, so
the leading eigenvectors are indeed sparse. This leads to the following
sparsity condition.

\renewcommand{\theconditionsp}{SP}
\begin{conditionsp}\label{condSP}
As $n\to\infty$, the radius $s_j$ of the weak-$\ell_r$ ball
satisfies $s_j\geq1$ and
\[
s_j^r \biggl[\frac{\log(\pn)}{n\lambda_j^4}
\biggr]^{\hf-r/4} = o\bigl(1\wedge\lambda_1^4\bigr)\qquad
\mbox{for $j=1,\ldots,\nspike$}.
\]
\end{conditionsp}

This type of condition also appeared in a previous
study of
individual eigenvector estimation in the multiple component spiked
covariance model~\cite{pajo07}.
The condition is, for example, satisfied if Condition~\ref{condGR} holds and the
largest spike $\lambda_1^2$ is bounded away from zero while the radii
$s_j$'s are all bounded above by an arbitrarily large constant. That
is, if there exists a constant $C>0$, such that $\lambda_1^2 \geq1/C$
and $s_j\leq C$ for all $j\leq\nspike$ and all $n$.

It is straightforward to verify that Conditions~\ref{condGR} and~\ref{condSP} are satisfied
by the special case in Section~\ref{secspecial}. We conclude this
part with an example.

\begin{example*}
When each $x_i$ collects noisy measurements of an underlying random
function on a regular 
grid, model (\ref{eqmodel}) becomes discretization of a functional
PCA model~\cite{fda}, and the $q_j$'s are discretized eigenfunctions.
When the eigenfunctions are smooth or have isolated singularities
either in themselves or in their derivatives, their wavelet
coefficients belong to some weak $\ell_r$ ball~\cite{dono93}. So do
the discrete wavelet transform of the $q_j$'s. Moreover, the radii of
the weak $\ell_r$ balls are determined by the underlying
eigenfunctions and are thus uniformly bounded as the size of the grid
$p$ gets larger. In this case, Condition~\ref{condSP} is satisfied when
Condition~\ref{condGR} holds and $\lambda_1^2$ is bounded away from zero. So, for
functional data of this type, we could always first transform to the
wavelet domain and then apply Algorithm~\ref{algspcait}.
\end{example*}

\subsection{Key quantities}
\label{secprelim}

We now introduce a few key quantities which appear later in the general
theoretical results.

The first quantity gives the rate at which we distinguish high from low
signal coordinates. Recall that $h(x) = x^2/(x+1)$. For $j=1,\ldots,
\nspike$, define
%
%
\begin{equation}
\label{eqtau-n} \tau_{nj} = \sqrt{\frac{\log(\pn)}{n h(\lambda_j^2)}}.
\end{equation}
According to~\cite{paul05}, up to a logarithmic factor, $\tau_{nj}^2$
can be interpreted as the average error per coordinate in estimating an
eigenvector with eigenvalue $\lambda_j^2+1$. Thus, a~coordinate can be
regarded as of high signal if at least one of the leading eigenvectors
is of larger magnitude on this coordinate compared to $\tau_{nj}$.
Otherwise, we call it a low signal coordinate.
We define $H(\beta)$ to be the set of high signal coordinates
%
%
\begin{equation}
\label{eqset-H} H = H(\beta) = \bigl\{ \nu\dvtx|q_{\nu j}|\geq\beta
\tau_{nj}\mbox{, for some $1\leq j\leq{\nspike}$} \bigr\}.
\end{equation}
Here, $\beta$ is a constant not depending on $n$, the actual value of
which will be specified in Theorem~\ref{thmerror}.
If $\nspike= 1$ and $q_1$ has $k$ nonzero entries all equal to
$1/\sqrt{k}$, then $H$ contains exactly these $k$ coordinates when $k
< n h(\lambda_1^2)/[\beta^2\log(\pn)]$,
which is guaranteed under Condition~\ref{condSP}.
In addition, let $L = \{1,\ldots,p\}\setminus H$ be the complement of
$H$. Here, $H$ stands for ``high,'' and $L$ for ``low''
(also recall $B$ in Algorithm~\ref{algdtspca}, where $B$ stands for ``big'').
The dependence of $H$, $L$ and $B$ on $n$ is suppressed for notational
convenience.

To understand the convergence rate of the subspace estimator stated
later in (\ref{eqerr-rate}), it is important to have an upper bound
for $\operatorname{card}(H)$, the cardinality of $H$.
To this end, define
%
%
\begin{equation}
\label{eqM-n} M_n = p \wedge\sum_{j=1}^{\nspike}
\frac{s_j^r}{\tau_{nj}^{r}}.
\end{equation}
The following lemma shows that a constant multiple of $M_n$ bounds
$\operatorname{card}(H)$.
The proof of the lemma is given in~\cite{supp}.
Thus, in the general result, $M_n$ plays the same role as the term
$\nspike s^r [nh(\lambda^2)/\log{p}]^{r/2}$ has played in Theorem
\ref{thmspecial}.
%
%
\begin{lemma}
\label{lemmaH}
For sufficiently large $n$, the cardinality of $H = H(\beta)$ satisfies
$\nspike\leq\operatorname{card}(H)\leq C M_n$ for a constant $C$
depending on
$\beta$ and $r$.
\end{lemma}
%

The last quantity we introduce is related to the ``parametric'' term in
the convergence rate. Let $\lambda^2_{\nspike+1} = 0$. For $j=1,\ldots,
\nspike$, define
%
%
\begin{equation}
\label{eqeps-nj} \eps_{nj}^2 = \frac{(\lambda_1^2+1)(\lambda
_{j+1}^2+1)}{(\lambda
_j^2-\lambda_{j+1}^2)^2}
\frac{\log(\pn)}{n}.
\end{equation}
So the second term of the upper bound in Theorem~\ref{thmspecial} is
$C_2\eps_{nm}^2$. For the interpretation of this quantity, we refer to
the discussion after Theorem~\ref{thmspecial}.

\subsection{Main results}
\label{secresult}
We turn to the statement of main theoretical results.

A key condition for the results is the asymptotic distinguishability
(AD) condition introduced below. Recall that all the spikes $\lambda
_j^2$ (hence\vspace*{1pt} all the leading eigenvalues) are allowed to depend on $n$.
The condition AD will guarantee that the largest few eigenvalues are
asymptotically well separated from the rest of the spectrum, and so the
corresponding principal subspace is distinguishable.

\begin{definition*}
We say that \textit{condition $\operatorname{AD}(j,\kappa)$ is
satisfied with
constant~$\kappa$}, if there exists a numeric constant $\kappa
\geq1$, such that for sufficiently large $n$, the gap between the
$j$th and the $(j+1)$th eigenvalues satisfies
\[
\lambda_j^2 - \lambda_{j+1}^2
\geq{\lambda_1^2}/{\kappa}. 
\]
\end{definition*}

We define $\operatorname{AD}(0,\kappa)$ and $\operatorname
{AD}(\nspike,\kappa)$ by letting
$\lambda_0^2 = \infty$, and $\lambda^2_{\nspike+1} = 0$.
So $\operatorname{AD}(0,\kappa)$ holds for any $\kappa\geq1$. %
Note that there is always some $1\leq j\leq\nspike$ such that
condition $\operatorname{AD}(j,\kappa)$ is satisfied. For instance,
$\operatorname{AD}(j,\kappa)$ is satisfied with some $\kappa$ for
the largest $j$ such that
$\lambda_j^2 \asymp\lambda_1^2$. When the spikes do not change with
$n$, condition $\operatorname{AD}(\nspike, \kappa)$ is satisfied
with any constant
$\kappa\geq\lambda_1^2/\lambda_{\nspike}^2$.
%

\subsubsection*{Rates of convergence for principal subspace estimation}
Recall definitions (\ref{eqalpha-n})--(\ref{eqstop}) and (\ref
{eqtau-n})--(\ref{eqeps-nj}).
The following theorem establishes the rate of convergence of the
principal subspace estimator obtained via Algorithm~\ref{algspcait}
under relaxed assumptions, which generalizes Theorem~\ref{thmspecial}.

%
%
\begin{theorem}
\label{thmerror}
Suppose Conditions~\ref{condGR} and~\ref{condSP} hold,
and condition $\operatorname{AD}(m,\kappa)$ is satisfied with some constant
$\kappa\geq1$ for the given subspace dimension~$m$.
Let the constants $\alpha,\gamma> 2\sqrt{3}$ in (\ref{eqalpha-n})
and (\ref{eqthr-nj}), and for $c = 0.9(\gamma-2\sqrt{3})$, let
$\beta= c/\sqrt{m}$ in $H$ (\ref{eqset-H}).
Then, there exist constants $C_0$, $C_1 = C_1(\thr, r, m, \kappa)$
and~$C_2$, such that for sufficiently large $n$, 
uniformly over $\class_n$,
with probability at least $1-C_0\pn^{-2}$, $\Ks\in[K, 2K]$ for
%
%
\begin{equation}
\label{eqK} {
K = \frac{\lambda_1^2+1}{\lambda_{m}^2- \lambda_{m+1}^2} \biggl[ \biggl
(1 +
\frac{1}{\log{2}} \biggr)\log{n} + 0\vee\log h\bigl(\lambda_1^2
\bigr) \biggr],}
\end{equation}
%
and the subspace estimator $\widehat{\subsp}_m^{({\Ks})} =
\operatorname{ran}(\widehat{Q}{}^{({\Ks})})$
satisfies
%
%
\begin{equation}
\label{eqerr-rate} 
L\bigl(\subsp_m, \widehat{
\subsp}_m^{({\Ks})}\bigr) \leq C_1 M_n
\tau_{nm}^2 + C_2\eps_{nm}^2
= o(1).
\end{equation}
%
\end{theorem}
%

Theorem~\ref{thmerror} states that for appropriately chosen threshold
levels and all thresholding function satisfying (\ref{eqeta}), after
enough iterations, Algorithm~\ref{algspcait} yields principal
subspace estimators whose errors are, with high probability, uniformly
bounded over $\class_n$ by a sequence of asymptotically vanishing
constants as $n\to\infty$. In addition, the probability that the
estimation error is not well controlled vanishes polynomially fast.
Therefore, the subspace estimators are uniformly consistent over
$\class_n$.\vadjust{\goodbreak}

The interpretation of the two terms in the error bound (\ref
{eqerr-rate}) is similar to those in Theorem~\ref{thmspecial}. Having
introduced those quantities in Section~\ref{secprelim}, we could
elaborate a little more on the first, that is, the ``nonparametric'' term.
By Theorem~\ref{thmexclude} below, when estimating $\subsp_m$,
Algorithm~\ref{algspcait} focuses only on the coordinates in $H$,
whose cardinality is $\operatorname{card}(H)= O(M_n)$.
Though $H$ does not appear explicitly in the rates, the rates
depend\vspace*{1pt} crucially on its cardinality which is further upper bounded by $M_n$.
Since $\tau_{nm}^2$ can be interpreted as the average error per
coordinate, the total estimation error accumulated over all coordinates
in $H$ is thus of order $O(M_n \tau_{nm}^2)$.
Moreover, as we will show later, the squared bias induced by focusing
only on $H$ is also of order $O(M_n\tau_{nm}^2)$. Thus, this term
indeed comes from the bias-variance tradeoff of the nonparametric
estimation procedure. The meaning of the second, that is, the
``parametric,'' term is the same as in Theorem~\ref{thmspecial}.
Finally, we note that both terms vanish as $n\to\infty$ under
Conditions~\ref{condGR},~\ref{condSP} and $\operatorname{AD}(m,\kappa)$.

The threshold levels $\alpha_n$ and $\thr_{nj}$ in (\ref
{eqalpha-n}) and (\ref{eqthr-nj}) as well as $\Ks$ in (\ref
{eqstop}) do not depend on unknown parameters.
So the estimation procedure achieves the rates adaptively over a wide
range of high-dimensional sparse settings.

In addition, (\ref{eqK}) implies that Algorithm~\ref{algspcait}
only needs a relatively small number of iterations to yield the desired
estimator. In particular, when the largest spike $\lambda_1^2$ is
bounded away from zero, (\ref{eqK}) shows that it suffices to have
$\Ks\asymp\log{n}$ iterations.
We remark that it is not critical to run precisely $K_s$
iterations. The result holds when we stop anywhere between $K$ and $2K$.

Theorem~\ref{thmerror} could also be extended to an upper bound for the
risk. Note that $\pn^{-2} = o(\tau_{nm}^2 \vee\eps_{nm}^2)$, and that
the loss function (\ref{eqloss}) is always bounded above by~$1$. The
following result is a direct consequence of Theorem~\ref{thmerror}.
%
%
\begin{corollary}
\label{corrisk}
Under the setup of Theorem~\ref{thmerror}, we have
\[
\sup_{\class_n} \Ex L\bigl(\subsp_m, \widehat{
\subsp}_m^{({\Ks})}\bigr) \leq{C_1}M_n
\tau_{nm}^2 + C_2 \eps_{nm}^2.
\]
\end{corollary}

\subsubsection*{Correct exclusion property}
We now switch to the model selection property of Algorithm~\ref{algspcait}.
By the discussion in Section~\ref{secalg}, an important motivation
for the iterative thresholding procedure is to trade bias for variance
by keeping low signal coordinates out of the orthogonal iterations.
More specifically,
it is desirable to restrict our effort to estimating those coordinates
in $H$ and simply estimating those coordinates in $L$ with zeros.

By construction, Algorithm~\ref{algdtspca} yields an initial matrix
with a lot of zeros, but Algorithm~\ref{algspcait} is at liberty to
introduce new nonzero coordinates. The following result shows that
with high probability all the nonzero coordinates introduced are in
the set $H$.

%
%
\begin{theorem}
\label{thmexclude}
Under the setup of Theorem~\ref{thmerror}, uniformly over $\class
_n$, with probability at least $1-C_0 \pn^{-2}$, for all $k=0,\ldots,
{\Ks}$, the orthonormal matrix $\widehat{Q}{}^{(k)}$ has zeros in
all its rows indexed by $L$, that is,
$\widehat{Q}{}^{(k)}_{L\cdot} = 0$.
\end{theorem}

We call the property in Theorem~\ref{thmexclude} ``correct
exclusion,'' because it ensures that all the low signal coordinates in
$L$ are correctly excluded from iterations.
In addition, Theorem~\ref{thmexclude} shows that the principal
subspace estimator is indeed spanned by a set of sparse loading
vectors, where all loadings in $L$ are exactly zero.

Note that the initial matrix $\widehat{Q}{}^{(0)}$ has all its
nonzero coordinates in $B$, which, with high probability, only selects
``big'' coefficients in the leading eigenvectors,
whose magnitudes are no less than $O([\log{\pn}/(n\lambda_m^4)]^{1/4})$.
On the other hand, the set $H$ includes all coordinates with magnitude
no less than $O([\log\pn/(n h(\lambda_m^2))]^{1/2})$. Thus, the
minimum signal strength for $H$ is of smaller order than that for $B$.
So, with high probability, $B$ is a subset of $H$ consisting only of
its coordinates with ``big'' signals.
Thus, though $\widehat{Q}{}^{(0)}$ excludes all the coordinates in
$L$, it only includes ``big'' coordinates in $H$ and fails to pick
those medium sized ones which are crucial for obtaining the convergence
rate (\ref{eqerr-rate}).
Algorithm~\ref{algspcait} helps to include more coordinates in $H$
along iterations and hence achieves (\ref{eqerr-rate}).

\subsubsection*{Rates of convergence for individual eigenvector estimation}
The primary focus of this paper is on estimating principal subspaces.
However, when an individual eigenvector, say $q_j$, is identifiable, it
is also of interest to see whether Algorithm~\ref{algspcait} can
estimate it well. The following result shows that for {$\Ks$ in (\ref
{eqstop}),}
the $j$th column of $\widehat{Q}{}^{({\Ks})}$ estimates $q_j$
well, provided that the $j$th eigenvalue is well separated from the
rest of the spectrum.

%
%
\begin{corollary}
\label{corevec-rate}
Under the setup of Theorem~\ref{thmerror}, suppose for some \mbox{$j\leq
m$}, both conditions $\operatorname{AD}(j-1,\kappa')$ and
$\operatorname{AD}(j,\kappa')$ are satisfied
for some constant $\kappa' < \lim_{n\to\infty}\lambda_1^2/(\lambda
_m^2-\lambda_{m+1}^2)$.
Then
uniformly over $\class_n$,
with probability at least $1-C_0\pn^{-2}$, $\widehat{q}_j^{({\Ks
})}$, the $j$th column of $\widehat{Q}{}^{({\Ks})}$, satisfies
\[
L\bigl(\operatorname{span}\{q_j\}, \operatorname{span}\bigl\{
\widehat{q}_j^{({\Ks})} \bigr\}\bigr) \leq C_1
M_n \tau_{nj}^2 + C_2 \bigl(
\eps_{n,j-1}^2\vee\eps_{nj}^2\bigr).
\]
Moreover, $\sup_{\class_n} \Ex L(\operatorname{span}\{q_j\},
\operatorname{span}\{\widehat{q}_j^{({\Ks})} \})$, the supremum risk over
$\class_n$, is also bounded by the right-hand side of the above inequality.
\end{corollary}

Corollary~\ref{corevec-rate} connects closely to the previous
investigation~\cite{pajo07} on estimating individual sparse leading
eigenvectors. Recall their loss function
$l(q_j, \tilde{q}_j) = \|q_j - \operatorname{sgn}(q_j'\tilde
{q}_j)\tilde{q}_j\|_{2}^2$.
Since
$\frac{1}{2}l(q_j, \tilde{q}_j) \leq
L(\operatorname{span}\{q_j\}, \operatorname{span}\{\tilde{q}_j\})
\leq\break
l(q_j, \tilde{q}_j)$, $l$ is
equivalent to the restriction of the loss function (\ref{eqloss}) to
one-dimensional subspaces.
Thus, Corollary~\ref{corevec-rate} implies that
\[
\sup_{\class_n} \Ex l\bigl(q_j, \widehat{q}_j^{({\Ks})}
\bigr)\leq C_1 M_n\tau_{nj}^2 +
C_2 \frac{(\lambda_1^2+1)(\lambda_j^2+1)}{
[(\lambda_{j-1}^2-\lambda_j^2)\wedge(\lambda_j^2-\lambda
_{j+1}^2)]^2 }\frac{ \log(\pn) }{n}.
\]

When the radii of the weak-$\ell_r$ balls grow at the same rate, that
is, $\max_j s_j\asymp\min_j s_j$, the upper bound in the last
display matches\vadjust{\goodbreak} the lower bound in Theorem 2 of~\cite{pajo07} up to a
logarithmic factor. Thus, when the $j$th eigenvalue is well separated
from the rest of the spectrum,
Algorithm~\ref{algspcait} yields a near optimal estimator of $q_j$ in
the adaptive rate minimax sense.

\subsection{Choice of $m$}
\label{secrank}
The main results in this section are stated with the assumption that
the subspace dimension $m$ is given.
In what follows, we discuss how to choose $m$ and also how to estimate
$\nspike$ based on data.

Recall $\ell^B_j$ defined after (\ref{eqthr-nj})
and the set $B$ in step 1 of Algorithm~\ref{algdtspca}.
Let 
%
%
\begin{equation}
\label{eqhat-nspike} \widehat{\nspike} = \max\bigl\{j\dvtx\ell^B_j
> 1 + \delta_{\operatorname{card}(B)} \bigr\}
\end{equation}
be an estimator for $\nspike$,
where for any positive integer $k$,
%
%
\begin{eqnarray}
\label{eqdelta-k} \delta_k & = & 2(\sqrt{k/ n} + t_k) +
(\sqrt{k/ n} + t_k )^2
\end{eqnarray}
with
\begin{eqnarray}
\label{eqt-k}
t_k^2 & = & \frac{6\log{\pn}}{n} +
\frac{2k(\log\pn+1)}{n}.
\end{eqnarray}
Then in Algorithm~\ref{algspcait}, for a large constant $\bar{\kappa
}$, we define
%
%
\begin{equation}
\label{eqm} m = \max\biggl\{j\dvtx1\leq j\leq\widehat{\nspike} \mbox{
and }
{\ell^B_1-1\over\ell^B_j - \ell^B_{j+1}} \leq\bar{\kappa} \biggr\}.
\end{equation}
Setting $\bar{\kappa}=15$ works well in simulation.
For a given dataset, such a choice of $m$ is intended to lead us to
estimate the largest principal subspace such that its eigenvalues
maintain a considerable gap from the rest of the spectrum.
Note that (\ref{eqm}) can be readily incorporated into Algorithm \ref
{algdtspca}:
we could compute the eigenvalues $\ell^B_j$ of $S_{BB}$ in step 2, and
then obtain $\widehat{\nspike}$ and $m$.

For $\widehat{\nspike}$ in (\ref{eqhat-nspike}), we have the
following results.
%
%
\begin{proposition}
\label{prophat-nspike}
Suppose Conditions~\ref{condGR} and~\ref{condSP} hold.
Let $\widehat{\nspike}$ be defined in (\ref{eqhat-nspike}) with $B$ obtained
by Algorithm~\ref{algdtspca} with $\alpha_n$ specified by (\ref
{eqalpha-n}) for some \mbox{$\alpha> 2\sqrt{3}$}.
Then, uniformly over $\class_n$, with probability at least $1 - C_0\pn
^{-2}$: 
\begin{longlist}[(3)]
\item[(1)] $\widehat{\nspike} \leq\nspike$;
\item[(2)] for any $m$ such that the condition
$\operatorname{AD}(m, \kappa)$ is satisfied with some constant
$\kappa$, $m\leq
\widehat{\nspike}$ when $n$ is sufficiently large;
\item[(3)] if the condition $\operatorname{AD}(\nspike,\kappa)$ is
satisfied with
some constant $\kappa$, $\widehat{\nspike} = \nspike$ when $n$ is
sufficiently large.
\end{longlist}
\end{proposition}
By claim (1), any $m\leq\widehat{\nspike}$ satisfies
$m\leq\nspike$ with high probability.
In addition, claim (2) shows that, for sufficiently large $n$, any $m$
such that $\operatorname{AD}(m,\kappa)$ holds
is no greater than $\widehat{\nspike}$.
Thus, when restricting to those $m\leq\widehat{\nspike}$, we do not
miss any $m$ such that Theorems~\ref{thmerror} and~\ref{thmexclude}
hold for estimating $\subsp_m$.
These two claims jointly ensure that we do not need to consider any
target dimension beyond $\widehat{\nspike}$.
Finally,
claim (3) shows that we recover\vadjust{\goodbreak} the exact number of spikes with high
probability for large samples when 
$\operatorname{AD}(\nspike, \kappa)$ is satisfied, that is, when
$\lambda_1^2
\asymp\lambda_{\nspike}^2$.
Note that this assumption was made in~\cite{pajo07}.

Turn\vspace*{1pt} to the justification of (\ref{eqm}).
We show later in Corollary~\ref{coreval} that for $1\leq j\leq
\nspike$, $(\ell^B_1-1)/(\ell^B_j - \ell^B_{j+1})$ estimates
$\lambda_1^2/(\lambda_j^2 - \lambda_{j+1}^2)$ consistently\vspace*{1pt} under
Conditions~\ref{condGR} and~\ref{condSP}.
[It is important that the condition $\operatorname{AD}(m,\kappa)$ is
not needed for
this result!] This implies that for $m$ in (\ref{eqm}), we have
$\lambda_1^2/(\lambda_m^2 - \lambda_{m+1}^2) \leq1.1\bar{\kappa}$
when $n$ is sufficiently large. Hence, the condition $\operatorname
{AD}(m,\kappa)$
is satisfied with the constant $\kappa= 1.1\bar{\kappa}$.
Therefore, the main theoretical results in Section~\ref{secresult}
remain valid when we set $m$ by (\ref{eqm}) in Algorithm~\ref{algspcait}.

\section{Computational complexity}
\label{seccomplex}
We now study the computational complexity of Algorithm \ref
{algspcait}. Throughout, we assume the same setup as in Section \ref
{secopt}, and restrict the calculation to the high probability event
on which the conclusions of Theorems~\ref{thmerror} and \ref
{thmexclude} hold.
For any matrix $A$, we use $\supp\{A\}$ to denote the index set of the
nonzero rows of $A$. 

Consider a single iteration, say, the $k$th.
In the multiplication step, the $(\nu, j)$th element of $T^{(k)}$,
$t_{\nu j}^{(k)}$, comes from the inner product of the
$\nu$th row of $S$ and the $j$th column of $\widehat{Q}{}^{(k-1)}$.
Though\vspace*{1pt} both are $p$-vectors, Theorem~\ref{thmexclude} asserts that
for any column of $\widehat{Q}{}^{(k-1)}$, at most $\operatorname
{card}(H)$ of
its entries are nonzero.
So if we know $\supp\{\widehat{Q}{}^{(k-1)} \}$, then $t_{\nu
j}^{(k)}$ can be calculated in $O(\operatorname{card}(H))$ flops,
and $T^{(k)}$ in $O(mp \operatorname{card}(H))$ flops.
Since $\supp\{\widehat{Q}{}^{(k-1)}\}$ can be obtained in $O
(mp)$ flops, the multiplication step can be completed in $O(mp
\operatorname{card}(H))$ flops.
Next, the thresholding step performs elementwise operation on
$T^{(k)}$, and hence can be completed in $O(mp)$ flops. Turn to the QR step.
First, we can obtain $\supp\{\widehat{T}{}^{(k)} \}$ in $O
(mp)$ flops. Then QR factorization can be performed on the reduced
matrix which only includes the rows in $\supp\{\widehat{T}{}^{(k)} \}$.
Since Theorem~\ref{thmexclude} implies $\supp\{\widehat
{T}{}^{(k)} \} = \supp\{\widehat{Q}{}^{(k)}\}\subset H$, the
complexity of this step is $O(m^2\operatorname{card}(H))$.
Since $m = O(p)$, the complexity of the multiplication step dominates,
and so the complexity of each iteration is $O(m p \operatorname{card}(H))$.
Theorem~\ref{thmerror} shows that {$\Ks$} iteration is enough.
Therefore, the overall complexity of Algorithm~\ref{algspcait} is
$O({\Ks}mp \operatorname{card}(H))$.

When the true eigenvectors are sparse, $\operatorname{card}(H)$ is of
manageable size.
In many realistic situations, $\lambda_1^2$ is bounded away from $0$
and so {$\Ks\asymp\log{n}$}. For these cases, Algorithm \ref
{algspcait} is scalable to very high dimensions.

We conclude the section with a brief discussion on parallel
implementation of Algorithm~\ref{algspcait}. In the $k$th iteration,
both matrix multiplication and elementwise thresholding can be computed
in parallel. For QR factorization, one needs only to communicate the
rows of $\widehat{T}{}^{(k)}$ with nonzero elements, the number
of which is no greater than $\operatorname{card}(H)$. Thus, the
overhead from
communication is $O(m \operatorname{card}(H))$ for each iteration,
and $O
({\Ks}m \operatorname{card}(H))$ in total. When the leading
eigenvectors are
sparse, $\operatorname{card}(H)$ is manageable, and parallel computing
of Algorithm
\ref{algspcait} is feasible.

\section{Numerical experiments}\vspace*{-6pt}
\label{secnumeric}
\subsection{Single spike settings}
\begin{figure}

\includegraphics{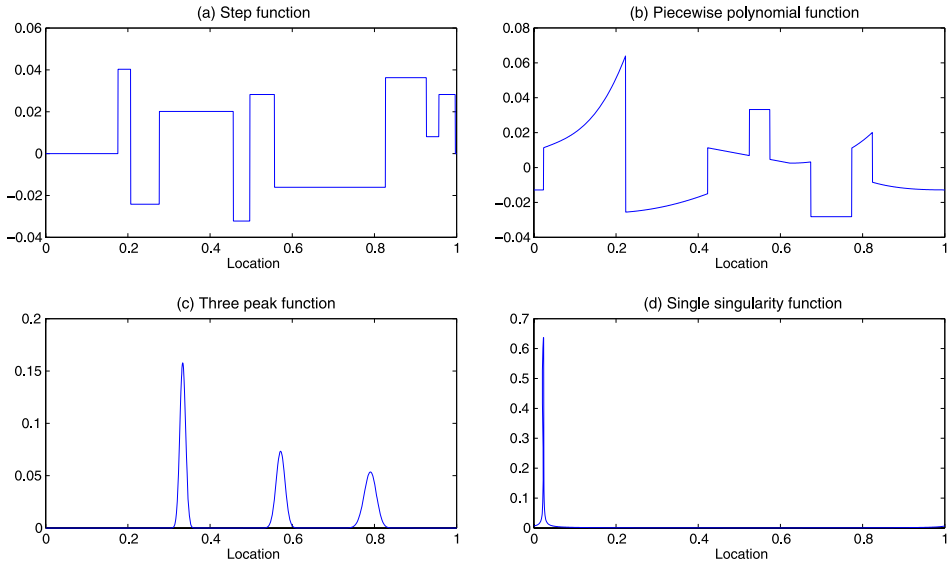}\vspace*{-5pt}

\caption{Four test vectors in the original domain: values at $p=2048$
equispaced points on $[0,1]$ of four test functions. \textup{(a)}
\texttt{step}: step
function, \textup{(b)} \texttt{poly}: piecewise polynomial function,
\textup{(c)} \texttt{peak}:
three-peak function and \textup{(d)} \texttt{sing}: single singularity
function.}\vspace*{-8pt}
\label{figtestorig}
\end{figure}

We first consider the case where each $x_i$ is generated by (\ref
{eqmodel}) with ${\nspike} = 1$.
Motivated by functional data with localized features, four test vectors
$q_1$ are considered, where
$q_1 = (f(1/p), \ldots, f(p/p))'$,
with $f$ one of the four functions in Figure~\ref{figtestorig}.
For each test vector, the dimension $p = 2048$, the sample size
$n=1024$ and $\lambda_{1}^2$ ranges in $\{100, 25, 10, 5, 2\}$.

%
%

\begin{figure}

\includegraphics{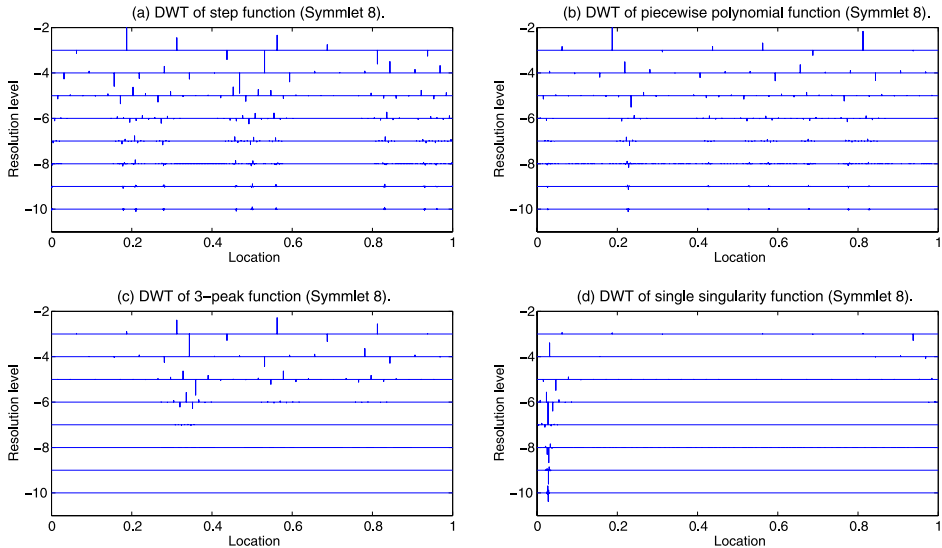}

\caption{Discrete wavelet transform of the four test vectors in
Figure \protect\ref{figtestorig}. In each plot, the length of each
stick is
proportional to the magnitude of the \texttt{Symmlet 8} wavelet
coefficient at the given location and resolution level.}
\label{figtestwave}
\end{figure}

Before applying any sparse PCA method, we transform the observed data
vectors into the wavelet domain using the Symmlet 8 basis \cite
{mall09}, and scale all the observations by $\widehat{\sigma}$ with
$\widehat\sigma^2$ given in (\ref{eqhat-sigma}). 
The multi-resolution plots of wavelet coefficients of the test vectors
are shown in Figure~\ref{figtestwave}. In the wavelet domain, the
four vectors exhibits different levels of sparsity, with \texttt{step} the
least sparse, and \texttt{sing} the most.

%
%

Table~\ref{tabsingle} compares the average loss of subspace
estimation over $100$ runs for each spike value and each test vector by
Algorithm~\ref{algspcait} (ITSPCA) with several existing methods:
augmented sparse PCA (AUGSPCA)
\cite{pajo07}, correlation augmented sparse PCA (CORSPCA)
\cite{nadl09} and diagonal thresholding sparse PCA (DTSPCA) given in
Algorithm~\ref{algdtspca}. 
For ITSPCA, we computed $\widehat{Q}{}^{(0)}$ by Algorithm \ref
{algdtspca}.
$\alpha_n$ and $\thr_{n1}$ are specified by (\ref{eqalpha-n}) and
(\ref{eqthr-nj}) with $\alpha=3$ and $\thr= 1.5$. These values are
smaller than those in theoretical results, but lead to better numerical
performance.
We stop iterating once $L(\operatorname{ran}(\widehat
{Q}{}^{(k)}),\operatorname{ran}(\widehat{Q}{}^{(k+1)})) \leq n^{-2}$.
Parameters in competing algorithms are all set to the values
recommended by their authors.\vadjust{\goodbreak}



%
%
\begin{table}
\def\arraystretch{0.9}
\caption{Comparison of sparse PCA methods in
single spike settings: average loss in estimation and size of selected
feature set}\label{tabsingle}
\begin{tabular*}{\tablewidth}{@{\extracolsep{\fill}}ld{3.0}crcrcrcr@{}}
\hline
&& \multicolumn{2}{c}{\textbf{ITSPCA}} & \multicolumn{2}{c}{\textbf{AUGSPCA}}
& \multicolumn{2}{c}{\textbf{CORSPCA}} & \multicolumn{2}{c@{}}{\textbf{DTSPCA}} \\[-4pt]
&& \multicolumn{2}{c}{\hrulefill} & \multicolumn{2}{c}{\hrulefill}
& \multicolumn{2}{c}{\hrulefill} & \multicolumn{2}{c@{}}{\hrulefill} \\
\textbf{Test vector} & \multicolumn{1}{c}{$\bolds{\lambda_1^2}$}
& \multicolumn{1}{c}{\textbf{Loss}} & \multicolumn{1}{c}{\textbf{Size}}
& \multicolumn{1}{c}{\textbf{Loss}} & \multicolumn{1}{c}{\textbf{Size}}
& \multicolumn{1}{c}{\textbf{Loss}} & \multicolumn{1}{c}{\textbf{Size}}
& \multicolumn{1}{c}{\textbf{Loss}} & \multicolumn{1}{c@{}}{\textbf{Size}}\\
\hline
\texttt{Step} & $100$ & $0.0061$ & $114.2$ & $0.0096$ & $96.5$ & $0.0055$ & $120.1$ &
$0.0275$ & $66.6$ \\
& $25$ & $0.0224$ & $76.3$ & $0.0362$ & $55.4$ & $0.0236$ & $73.9$ &
$0.0777$ & $38.3$ \\
& $10$ & $0.0470$ & $53.4$ & $0.0710$ & $37.4$ & $0.0551$ &
$45.9$ & $0.1494$ & $24.1$ \\
& $5$ & $0.0786$ & $45.5$ & $0.1370$ & $23.7$ & $0.1119$ & $28.7$ &
$0.2203$ & $17.1$ \\
& $2$ & $0.1921$ & $25.4$ & $0.3107$ & $11.4$ & $0.3846$ & $15.2$ &
$0.4518$ & $9.7$ \\
[6pt]
\texttt{Poly} & $100$ & $0.0060$ & $83.1$ & $0.0088$ & $66.5$ & $0.0051$ & $92.0$ &
$0.0191$ & $49.2$ \\
& $25$ & $0.0175$ & $52.4$ & $0.0254$ & $41.4$ & $0.0173$ & $53.1$ &
$0.0540$ & $28.7$ \\
& $10$ & $0.0346$ & $38.7$ & $0.0527$ & $27.5$ & $0.0404$ &
$34.0$ & $0.0959$ & $20.5$ \\
& $5$ & $0.0588$ & $30.7$ & $0.0844$ & $20.2$ & $0.0684$ & $24.6$ &
$0.1778$ & $14.0$ \\
& $2$ & $0.1317$ & $20.0$ & $0.2300$ & $10.3$ & $0.2155$ & $16.3$ &
$0.3370$ & $8.1$ \\
[6pt]
\texttt{Peak} & $100$ & $0.0019$ & $45.7$ & $0.0032$ & $39.6$ & $0.0016$ & $51.2$ &
$0.0075$ & $32.8$ \\
& $25$ & $0.0071$ & $34.1$ & $0.0099$ & $29.9$ & $0.0069$ & $35.2$ &
$0.0226$ & $24.3$ \\
& $10$ & $0.0158$ & $28.0$ & $0.0222$ & $23.8$ & $0.0165$ &
$27.3$ & $0.0592$ & $18.6$ \\
& $5$ & $0.0283$ & $24.7$ & $0.0449$ & $19.6$ & $0.0320$ & $22.5$ &
$0.1161$ & $14.1$ \\
& $2$ & $0.0927$ & $20.8$ & $0.1887$ & $9.9$ & $0.1176$ & $14.6$ &
$0.2702$ & $8.8$ \\
[6pt]
\texttt{Sing} & $100$ & $0.0016$ & $38.0$ & $0.0025$ & $33.2$ & $0.0014$ & $43.6$ &
$0.0070$ & $26.3$ \\
& $25$ & $0.0068$ & $27.1$ & $0.0095$ & $23.1$ & $0.0060$ & $31.8$ &
$0.0237$ & $17.5$ \\
& $10$ & $0.0161$ & $20.3$ & $0.0233$ & $16.6$ & $0.0154$ &
$20.9$ & $0.0377$ & $13.6$ \\
& $5$ & $0.0279$ & $17.3$ & $0.0372$ & $13.2$ & $0.0313$ & $15.2$ &
$0.0547$ & $12.7$ \\
& $2$ & $0.0631$ & $15.2$ & $0.0792$ & $10.9$ & $0.0652$ & $13.0$ &
$0.2025$ & $8.8$ \\
\hline
\end{tabular*}   \vspace*{-3pt}
\end{table}

From Table~\ref{tabsingle}, ITSPCA and CORSPCA outperform the
other two methods in all settings.
Between the two, CORSPCA only wins by small margins when the spike
values are large. Otherwise, ITSPCA wins, sometimes with large margins.
For the same algorithm at the same spike value, the sparser the signal,
the smaller the estimation error.

Table~\ref{tabsingle} also presents the average sizes of the sets of
selected coordinates. While all methods yield sparse PC loadings,
AUGSPCA and DTSPCA seem to select too few coordinates, and thus
introduce too much bias. ITSPCA and CORSPCA apparently result in a
better bias-variance tradeoff.

\subsection{Multiple spike settings} Next, we simulated data vectors
using model (\ref{eqmodel}) with ${\nspike} = 4$. The $q_j$ vectors
are taken to be the four test vectors used in single spike settings, in
the same order as in Figure~\ref{figtestorig}, up to
orthonormalization.\footnote{The four test vectors are shifted such
that the inner product of any pair is close to $0$. So the vectors
after orthonormalization are visually indistinguishable from those in
Figure~\ref{figtestorig}.} We tried four different configurations of
the spike values $(\lambda_1^2,\ldots,\lambda_4^2)$, as specified in
the first column of Table~\ref{tabmultiple}. For each configuration
of spike values, the dimension is $p=2048$, and the sample size is $n=1024$.

%
%
\begin{table}
\def\arraystretch{0.9}
\caption{Comparison of sparse PCA methods in
multiple spike settings: average loss in estimation}\label{tabmultiple}
\begin{tabular*}{\tablewidth}{@{\extracolsep{\fill}}lccccc@{}}
\hline
& & \multicolumn{4}{c@{}}{$\bolds{L(\subsp_m, \widehat{\subsp}_m)}$}\\[-4pt]
& & \multicolumn{4}{c@{}}{\hrulefill}\\
$\bolds{(\lambda_1^2,\lambda_2^2,\lambda_3^2,\lambda_4^2)}$ & $\bolds{m}$
& \textbf{ITSPCA} & \textbf{AUGSPCA} & \textbf{CORSPCA} & \textbf{DTSPCA} \\
\hline
$(100, 75, 50, 25)$ & $1$ & $0.0216$ & $0.0260$ & $0.0240$ & $0.0378$ \\
& $2$ & $0.0180$ & $0.0213$ & $0.0214$ & $0.0308$ \\
& $3$ & $0.0094$ & $0.0129$ & $0.0126$ & $0.0234$ \\
& $4$ & $0.0087$ & $0.0122$ & $0.0181$ & $0.0235$ \\
[6pt]
$(60, 55, 50, 45)$ & $1$ & $0.3100$ & $0.2588$ & $0.2548$ & $0.2831$ \\
& $2$ & $0.2675$ & $0.2045$ & $0.2095$ & $0.2349$ \\
& $3$ & $0.1844$ & $0.1878$ & $0.1872$ & $0.1968$ \\
& $4$ & $0.0157$ & $0.0203$ & $0.0178$ & $0.0333$ \\
[6pt]
$(30, 27, 25, 22)$ & $1$ & $0.3290$ & $0.2464$ & $0.2495$ & $0.2937$ \\
& $2$ & $0.3147$ & $0.2655$ & $0.2882$ & $0.3218$ \\
& $3$ & $0.1740$ & $0.1662$ & $0.1708$ & $0.1821$ \\
& $4$ & $0.0270$ & $0.0342$ & $0.0338$ & $0.0573$ \\
[6pt]
$(30, 20, 10, 5)$ & $1$ & $0.0268$ & $0.0392$ & $0.0380$ & $0.0658$ \\
& $2$ & $0.0237$ & $0.0353$ & $0.0391$ & $0.0605$ \\
& $3$ & $0.0223$ & $0.0336$ & $0.0372$ & $0.0599$ \\
& $4$ & $0.0298$ & $0.0414$ & $0.0717$ & $0.0638$ \\
\hline
\end{tabular*} \vspace*{-3pt}
\end{table}

For each simulated dataset, we estimate $\subsp_m$ for $m=1,2,3$ and
$4$. The last four columns of Table~\ref{tabmultiple} present the
losses in estimating subspaces, averaged over $100$ runs, using the
same sparse PCA methods as in single spike settings. For ITSPCA, we set
the thresholds $\{\thr_{nj},j=1,\ldots, 4\}$ as in (\ref{eqthr-nj})
with $\thr= 1.5$.
All other implementation details are the same. Again, we used
recommended values for parameters in all other competing methods.

The simulation results reveal two interesting phenomena. First, when
the spikes are relatively well separated (the first and the last blocks
of Table~\ref{tabmultiple}), all methods yield decent estimators of
$\subsp_m$ for all values of $m$, which implies that the individual
eigenvectors are also estimated well.
In this case, ITSPCA always outperforms the other three competing
methods. Second, when the spikes are not so well separated (the middle
two blocks, with $m=1,2$ or $3$), no method leads to decent subspace estimator.
However, all methods give reasonable estimators for $\subsp_4$ because
$\lambda_4^2$ in both cases are well above $0$. This implies that,
under such settings, we fail to recover individual eigenvectors, but we
can still estimate $\subsp_4$ well. ITSPCA again gives the smallest
average losses.\vadjust{\goodbreak}
In all configurations, the estimated number of spikes $\widehat
{\nspike}$ in (\ref{eqhat-nspike}) and the data-based choice of
$m$ in (\ref{eqm}) with $\bar{\kappa} = 15$ consistently picked
$m=\widehat{\nspike} = 4$ in all simulated datasets.
Therefore, we are always led to estimating the ``right'' subspace
$\subsp_4$, and ITSPCA performs favorably over the competing
methods.

In summary, simulations under multiple spike settings not only
demonstrate the competitiveness of Algorithm~\ref{algspcait}, but
also suggest:
\begin{longlist}[(2)]
\item[(1)] The quality of principal subspace estimation depends on the
gap between successive eigenvalues, in addition to the sparsity of eigenvectors;
\item[(2)] Focusing on individual eigenvectors can be misleading for
the purpose of finding low-dimensional projections.\vspace*{-2pt}
\end{longlist}

\section{Proof}
\label{secproof}
This section is devoted to the proofs of Theorems~\ref{thmerror} and
\ref{thmexclude}.
We state the main ideas in Section~\ref{secoutline-proof} and divide
the proof into three major steps, which are then completed in sequel in
Sections~\ref{secoracle-cov}--\ref{secora-act-id}.
Others results in Section~\ref{secresult} are proved in the
supplementary material~\cite{supp}.\vspace*{-2pt}

\subsection{Main ideas and outline of proof}
\label{secoutline-proof}

The proof is based on an \textit{oracle sequence approach}, the main
ideas of which are as follows.
First, assuming oracle knowledge of the 
set $H$, we construct a sequence of $p\times m$ orthonormal matrices $\{
\widehat{Q}{}^{(k),\orac},k\geq0\}$.
Then 
we study how fast the sequence converges, and how well each associated
column subspace approximates the principal subspace $\subsp_m$ of interest.
Finally,\vspace*{1pt} we show that, with high probability, the first $\Ks$ terms of
the oracle sequence is exactly the sequence $\{\widehat{Q}{}^{(k)},
0\leq k\leq\Ks\}$ obtained by Algorithm~\ref{algspcait}.
The actual estimating sequence thus inherits from the oracle sequence
various properties in terms\vadjust{\goodbreak} of estimation error and number of steps
needed to achieve the desired error rate.
The actual sequence mimics the oracle because the thresholding step
forces it to only consider the high signal coordinates in $H$.


In what follows, we first construct the oracle sequence and then lay
out a road map of the proof. Here and after, we use an extra
superscript ``$\orac$'' to indicate oracle quantities. For example,
$\widehat{Q}{}^{(k),\orac}$ denotes the $k$th orthonormal matrix
in the oracle sequence.\vspace*{-2pt}



\subsubsection*{Construction of the oracle sequence}

First, we construct $\widehat{Q}{}^{(0),\orac}$ using an oracle
version of Algorithm~\ref{algdtspca}, where the set $B$ is replaced
by its oracle version $B^\orac= B\cap H$. This ensures that $\widehat
{Q}{}^{(0),\orac}_{L\cdot} = 0$.

To construct the rest of the sequence,
suppose that the $p$ features are organized (after reordering) in such
a way that those in $H$ always have smaller indices than those in $L$,
and that within $H$, those in $B^\orac$ precede those not.
Define the oracle sample covariance matrix
%
%
\begin{equation}
\label{eqoracle-S} S^{\orac} = \lleft[\matrix{ S_{HH} & 0
\cr
0
& I_{LL}} \rright].
\end{equation}
Here, $I_{LL}$ is the identity matrix of dimension $\operatorname{card}(L)$.
Then,
the matrices $\{\widehat{Q}{}^{(k),\orac}, k\geq0\}$ are
obtained via an oracle version of Algorithm~\ref{algspcait}, in which
the initial matrix is $\widehat{Q}{}^{(0),\orac}$, and $S$ is
replaced by $S^{\orac}$.\vspace*{-2pt}

%
%
\begin{remark}
This formal construction does not guarantee that $\widehat
{Q}{}^{(k),\orac}$ has full column rank or that $\widehat
{Q}{}^{(k),\orac
}_{L\cdot} = 0$ for all $k$. Later, Lemma \ref
{lemmainitial}, Proposition~\ref{propevolve} and Lemma \ref
{lemmaexclude} show that these statements are true with high
probability for all {$k\leq\Ks$}.\vspace*{-2pt}
\end{remark}

\subsubsection*{Major steps of the proof}
In the $k$th iteration of the oracle Algorithm~\ref{algspcait},
denote the matrices obtained after multiplication and thresholding by
%
%
\begin{eqnarray}
\label{eqoracle-T}
T^{(k),\orac} & = &
S^{\orac} \widehat{Q}{}^{(k-1),\orac} = \bigl(t_{\nu j}^{(k),\orac}
\bigr)\quad \mbox{and}
\nonumber\\[-9.5pt]\\[-9.5pt]
\widehat{T}{}^{(k),\orac} & = & \bigl(\widehat{t}_{\nu j}^{(k),\orac}
\bigr) \qquad\mbox{with $\widehat{t}_{\nu j}^{(k),\orac} = \eta
\bigl({t}_{\nu j}^{(k),\orac}, \thr_{nj}\bigr)$.}\nonumber
\end{eqnarray}
Further denote the QR factorization of $\widehat{T}{}^{(k),\orac
}$ by $\widehat{T}{}^{(k),\orac} = \widehat{Q}{}^{(k),\orac
}\widehat{R}{}^{(k),\orac}$. Last but not least, let $\widehat
{\subsp}_m^{(k),\orac} = \operatorname{ran}(\widehat{Q}{}^{(k),\orac})$.

A joint proof of Theorems~\ref{thmerror} and~\ref{thmexclude} can
then be completed by the following three major steps:
\begin{longlist}[(3)]
\item[(1)] show that the principal subspace of $S^{\orac}$ with
dimension $m$, denoted by $\widehat{\subsp}_m^{\orac}$, satisfies
the error bound in (\ref{eqerr-rate}) for estimating $\subsp_m$;
\item[(2)] show that for $K$ in (\ref{eqK}), $\Ks\in[K, 2K]$ and
that the approximation error of $\widehat{\subsp}_m^{(k),\orac
}$ to $\widehat{\subsp}_m^{\orac}$ for all $k\geq K$ also satisfies
the bound in (\ref{eqerr-rate});
\item[(3)] show that $\widehat{Q}_{L\cdot}^{(k),\orac} = 0$
for all $k\leq2K$, and that the oracle and the actual estimating
sequences are identical up to $2K$ iterations.\vadjust{\goodbreak}
\end{longlist}
%
In each step, we only need the result to hold with high probability.
By the triangle inequality, steps 1 and 2 imply that the error of
$\widehat{\subsp}_m^{{(\Ks)},\orac}$ in\vspace*{1pt} estimating $\subsp
_m$ satisfies~(\ref{eqerr-rate}). Step 3 shows this is also the case for the actual
estimator $\widehat{\subsp}_m^{{(\Ks)}}$.
It also implies the correct exclusion property in Theorem~\ref{thmexclude}.

In what follows, we complete the three steps in Sections \ref
{secoracle-cov}--\ref{secora-act-id}.

\subsection{Principal subspace of $S^{\orac}$}
\label{secoracle-cov}

To study how well the principal subspace of $S^{\orac}$ approximates
$\subsp_m$,
we break into a ``bias'' part and a ``variance'' part.

Consider the ``bias'' part first.
Define the oracle covariance matrix
%
%
\begin{equation}
\label{eqoracle-Sigma} \Sigma^{\orac} = \lleft[\matrix{
\Sigma_{HH} & 0
\cr
0 & I_{LL}} \rright],
\end{equation}
which is the expected value of $S^{\orac}$. The following lemma gives
the error of the principal subspace of $\Sigma^{\orac}$ in
approximating $\subsp_m$, which could be regarded as the ``squared
bias'' induced by feature selection.

%
%
\begin{lemma}
\label{lemmaoracle-Sigma}
Let the eigenvalues of $\Sigma^{\orac}$ be
$\eval_{1}^{\orac}\geq\cdots\geq\eval_{\nspike}^{\orac}\geq\cdots\geq0$
and $\{
q_1^{\orac},\ldots,q_{\nspike}^{\orac}\}$ be a set of first $\nspike
$ eigenvectors. Denote $Q^{\orac} = [q_1^{\orac},\ldots, q_m^{\orac}]$.
Then, uniformly over $\class_n$:
\begin{longlist}[(2)]
\item[(1)] $|\eval_{j}^{\orac} - (\lambda_j^2+1)| / \lambda_1^2
\to0$ as $n\to\infty$, for $j=1,\ldots,\nspike+1$, with $\lambda
_{\nspike+1}^2=0$;
%
\item[(2)] for sufficiently large $n$, 
$Q^{\orac}_{L\cdot} = 0$ and for ${\subsp}^{\orac}_m =
\operatorname{ran}(Q^{\orac})$, there exists a constant $C = C(m, r,
\kappa)$, s.t.
$L(\subsp_m, \subsp^{\orac}_m) \leq C M_n \tau_{nm}^2$.
\end{longlist}
\end{lemma}

A proof is given in the supplementary material~\cite{supp}. Weyl's theorem
(\cite{stsu90}, Corollary 4.4.10) and Davis--Kahn's $\sin\theta$ theorem
\cite{daka70} are the key ingredients in the proof here, and also in
the proofs of Lemmas~\ref{lemmaoracle-S} and~\ref{lemmainitial}.
Here, claim (1) only requires Conditions~\ref{condGR} and~\ref{condSP}, but not the
condition $\operatorname{AD}(m,\kappa)$.


Turn to the ``variance'' part. We check how well the principal subspace
of $S^{\orac}$ estimates $\subsp^{\orac}$. 
Since $\Sigma^{\orac} = \Ex[S^{\orac}]$, the error here is
analogous to ``variance.''

%
%
\begin{lemma}
\label{lemmaoracle-S}
Let the eigenvalues of $S^{\orac}$ be $\widehat\eval_{1}^{\orac
}\geq\cdots\geq\widehat\eval_{\nspike}^{\orac}\geq\cdots\geq
0$ and $\{\widehat{q}_1^{\orac},\ldots,\widehat{q}_{\nspike}^{\orac
}\}$ be a set of first $\nspike$ eigenvectors. Denote $\widehat
{Q}{}^{\orac} = [\widehat{q}_1^{\orac},\ldots, \widehat{q}_m^{\orac}]$.
Then, uniformly over $\class_n$, with probability at least $1- C_0\pn^{-2}$:
\begin{longlist}[(2)]
\item[(1)] $|\widehat\eval_j^{\orac} - \eval_j^{\orac}|/\lambda
_1^2 \to0$ as $n\to\infty$, for $j=1,\ldots, \nspike+1$;
\item[(2)] for sufficiently large $n$, 
$\widehat{Q}{}^{\orac}_{L\cdot} = 0$, and for $\widehat{\subsp
}{}^{\orac}_m = \operatorname{ran}(\widehat{Q}{}^{\orac})$,
there exist constants $C_1 = C_1(m,r,\kappa)$ and $C_2$, s.t.
\[
L\bigl({\subsp}^{\orac}_m, \widehat{
\subsp}{}^{\orac}_m\bigr) \leq C_1 M_n
\tau_{nm}^2/\log(\pn) + C_2
\eps_{nm}^2.
\]
\end{longlist}
\end{lemma}

A proof is given in the supplementary material~\cite{supp}. Again,
claim (1) does
not require the condition $\operatorname{AD}(m,\kappa)$.
By the triangle inequality, the above two lemmas imply the error in
estimating $\subsp_m$ with $\widehat{\subsp}_m^{\orac}$ satisfies
the bound in (\ref{eqerr-rate}). 


\subsection{Properties of the oracle sequence}
\label{secoracle}

In step 2, we study properties of the oracle sequence. For $K$ in
(\ref{eqK}), the goal is to show that, with high probability, for all
$k\geq K$, the error of the oracle subspace estimator
$\widehat{\subsp}_m^{(k),\orac}$ in approximating $\widehat
{\subsp}_m^{\orac}$ satisfies in (\ref{eqerr-rate}).
To this end, characterization of the oracle sequence evolution in
Proposition~\ref{propevolve} below plays the key role.

\subsubsection*{The initial point}
We start with the initial point. 
Let
%
%
\begin{equation}
\label{eqrho} \rho= {\widehat{\eval}{}^{\orac}_{m+1}} / {\widehat{
\eval}{}^{\orac}_{m}}
\end{equation}
denote the ratio between the $(m+1)$th and the $m$th largest
eigenvalues of~$S^{\orac}$.
The following lemma shows that $\widehat{Q}{}^{(0),\orac}$ is
orthonormal and is a good initial point for (oracle) Algorithm \ref
{algspcait}.

%
%
\begin{lemma}
\label{lemmainitial}
Uniformly over $\class_n$, with probability at least $1-C_0\pn^{-2}$:
\begin{longlist}[(4)]
\item[(1)] $B^{\orac} = B$;
\item[(2)] $| \eval_j(S_{B^{\orac}B^{\orac}})\vee1 -
\widehat\eval_j^{\orac} | / \lambda_1^2 \to0$ as $n\to\infty$,
for $j=1,\ldots, \nspike+1$;
\item[(3)] 
for sufficiently large $n$, 
$\widehat{Q}{}^{(0),\orac}$ has full column rank, and $L(\widehat
{\subsp}{}^{\orac}_m, \widehat{\subsp}_m^{(0),\orac}) \leq
(1-\rho)^2/5$;
\item[(4)] for sufficiently large $n$, $\Ks\in[K, 2K]$.
\end{longlist}
\end{lemma}
A proof is given in the supplementary material~\cite{supp}.
Here, claims (1) and (2) do not require the condition $\operatorname
{AD}(m,\kappa)$.
In claim (3), the bound $(1-\rho)^2/5$ is much larger than that in
(\ref{eqerr-rate}). For instance, if $\lambda_m^2+1 \asymp\lambda
_m^2 - \lambda_{m+1}^2$, Lemmas~\ref{lemmaoracle-Sigma} and~\ref
{lemmaoracle-S} imply that $(1-\rho^2)/5 \asymp1$ with high probability.

Claims (1) and (2) here, together with claims (1) of Lemmas \ref
{lemmaoracle-Sigma} and~\ref{lemmaoracle-S}, lead to the following
result on consistent estimation of $\lambda_1^2/(\lambda_j^2-\lambda
_{j+1}^2)$ and $\lambda_j^2$, the proof of which is given in the
supplementary material~\cite{supp}.

%
%
\begin{corollary}
\label{coreval}
Suppose Conditions~\ref{condGR} and~\ref{condSP} hold, and let $\ell^B_j = \eval
_j(S_{BB})\vee1$.
For $1\leq j\leq\nspike$, if $\lim_{n\to\infty} \lambda
_1^2/(\lambda_j^2-\lambda_{j+1}^2) < \infty$, then
\[
\lim_{n\to\infty}\frac{(\ell^B_1-1)/(\ell^B_j-\ell^B
_{j+1})}{\lambda_1^2/(\lambda_j^2-\lambda_{j+1}^2)} = 1\qquad \mbox{a.s.}
\]
Otherwise, $\lim_{n\to\infty} (\ell^B_1-1)/(\ell^B_j-\ell^B
_{j+1}) = \lim_{n\to\infty} \lambda_1^2/(\lambda_j^2-\lambda
_{j+1}^2)= \infty$, a.s.

If further the condition $\operatorname{AD}(m,\kappa)$ holds for some
$m\leq
\nspike$ and $\kappa> 0$, then
$\lim_{n\to\infty} (\ell^B_j-1)/\lambda_j^2 = 1$, a.s.,
for $1\leq j\leq m$.
\end{corollary}

%


\subsubsection*{Evolution of the oracle sequence}
Next, we study the evolution of the oracle sequence. 
Let $\cangle^{(k)}\in[0, \pi/2]$ be the largest canonical\vadjust{\goodbreak}
angle between the subspaces $\widehat{\subsp}{}^{\orac}_m$ and
$\widehat{\subsp}_m^{(k),\orac}$. By the discussion after
(\ref{eqloss}), we have
%
%
\begin{equation}
\label{eqangle-loss} \sin^2 \cangle^{(k)} = L \bigl(
\widehat{\subsp}_m^{\orac}, \widehat{\subsp}_m^{(k),\orac}
\bigr).
\end{equation}
%
The following proposition describes the evolution of $\cangle^{(k)}$
over iterations.

%
%
\begin{proposition}
\label{propevolve}
Let $n$ be sufficiently large. On the event such that the conclusions
of Lemmas~\ref{lemmaoracle-Sigma}--\ref{lemmainitial} hold,
uniformly over $\class_n$, for all $k\geq1$:
%
\begin{longlist}[(2)]
\item[(1)] $\widehat{Q}{}^{(k),\orac}$ is orthonormal, and
$\cangle^{(k)}$ satisfies
%
%
\begin{equation}
\label{eqevolve} \sin\cangle^{(k)} \leq\rho\tan\cangle^{(k-1)} +
\omega\sec\cangle^{(k-1)},
\end{equation}
where $\omega= (\widehat{\eval}_m^{\orac})^{-1} [\operatorname
{card}(H)\sum_{j=1}^m \thr_{nj}^2 ]^\hf$;

\item[(2)] for any $a\in(0,1/2]$, if
%
%
\begin{equation}
\label{eqsubsp-nb} \sin^2\cangle^{(k-1)}
\leq{1.01(1-a)^{-2}\omega^2} {(1-\rho)^{-2}},
\end{equation}
then so is $\sin^2{\cangle^{(k)}}$. Otherwise,
%
%
\begin{equation}
\label{eqsubsp-rate} \sin^2\cangle^{(k)}/\sin^2
\cangle^{(k-1)} \leq\bigl[1-a(1-\rho)\bigr]^2.
\end{equation}
\end{longlist}
\end{proposition}

A proof is given in the supplementary material~\cite{supp}, the key
ingredient of
which is Wedin's $\sin\theta$ theorem for singular subspaces~\cite{wedi72}.
The recursive inequality (\ref{eqevolve}) characterizes the evolution
of the angles $\cangle^{(k)}$, and hence of the oracle
subspace~$\widehat{\subsp}_m^{(k),\orac}$. It is the foundation of
claim (2) in the current proposition and of Proposition~\ref
{proporacle-error} below.

By (\ref{eqangle-loss}), inequality (\ref{eqsubsp-rate}) 
gives the rate at which the approximation error $L(\widehat{\subsp
}_m^{\orac}, \widehat{\subsp}_m^{(k),\orac})$ decreases.
For a given $a\in(0,1/2]$, the rate is maintained until the error
becomes smaller than $1.01(1-a)^{-2}\omega^2(1-\rho)^{-2}$. Then the
error continues to decrease, but at a slower rate, say, with $a$
replaced by $a/2$ in (\ref{eqsubsp-rate}), until (\ref{eqsubsp-nb})
is satisfied with $a$ replaced by $a/2$. The decrease continues at
slower and slower rate in this fashion until the approximation error
falls into the interval $[0, {1.01\omega^2}/{(1-\rho)^2}]$, and
remains inside thereafter.

Together with Lemma~\ref{lemmainitial}, Proposition \ref
{propevolve} also justifies the previous claim that elements of the
oracle sequence are orthonormal with high probability.


\subsubsection*{Convergence}

Finally,
we study how fast the oracle sequence converges to a stable subspace
estimator, and how good this estimator is.

To define convergence of the subspace sequence $\{\widehat{\subsp
}_m^{(k),\orac}, k\geq0 \}$, we first note that $1.01\omega
^2/(1-\rho)^{2}$ is almost the smallest possible value of $L(\widehat
{\subsp}_m^{\orac}, \widehat{\subsp}_m^{(k),\orac})$ that
(\ref{eqevolve}) could imply. Indeed, when $\sin\cangle^{(k)}$
converges and is small, we have $\sin\cangle^{(k)}\approx
\sin\cangle^{(k-1)}$, and $\cos\cangle^{(k)}\approx1$.
Consequently, (\ref{eqevolve}) reduces to
\[
\sin\cangle^{(k)}\leq\bigl(\rho\sin\cangle^{(k)}+\omega\bigr)
\bigl(1+o(1)\bigr).
\]
So, $L(\widehat{\subsp}_m^{\orac}, \widehat{\subsp}_m^{(k),\orac
}) = \sin^2\cangle^{(k)}\leq(1+o(1)) \omega
^2/(1-\rho)^{2}$.
In addition, Lemma~\ref{lemmaoracle-S}
suggests that we can stop the iteration\vadjust{\goodbreak} as soon as $L(\widehat{\subsp
}_m^{\orac}, \widehat{\subsp}_m^{(k),\orac})$ becomes smaller
than a constant multiple of $\eps_{nm}^2$, for we always get an error
of order $O(\eps_{nm}^2)$ for estimating $\subsp_m$, even if we
use $\widehat{\subsp}_m^{\orac}$ directly.
In observation of both aspects, we say that $\widehat{\subsp
}_m^{(k),\orac}$ has \textit{converged} if
%
%
\begin{equation}
\label{eqoracle-conv} L\bigl(\widehat{\subsp}_m^{\orac},
\widehat{\subsp}_m^{(k),\orac
}\bigr) \leq\max\biggl\{
\frac{1.01}{(1-n^{-1})^2}\frac{\omega^2}{(1-\rho)^2}, \eps_{nm}^2 \biggr
\}.
\end{equation}
%
On the event that conclusions of Lemmas \ref
{lemmaoracle-Sigma}--\ref{lemmainitial} hold, we have $\omega
^2/(1-\rho)^2 = O( M_n\tau_{nm}^2)$.
Under definition (\ref{eqoracle-conv}), for $K$ in {(\ref{eqK})},
the following proposition shows that it takes $K$ iterations for the
oracle sequence to converge, and for all $k\geq K$, the error of
approximating $\widehat{\subsp}_m^{\orac}$ by $\widehat{\subsp
}_m^{(k),\orac}$ satisfies (\ref{eqerr-rate}).

%
%
\begin{proposition}
\label{proporacle-error}
For sufficiently large $n$, 
on the event such that the conclusions of Lemmas \ref
{lemmaoracle-Sigma}--\ref{lemmainitial} hold, uniformly over $\class
_n$, it takes at most
$K$ steps for the oracle sequence to converge. In addition, there exist
constants $C_1 = C_1(\thr,r,m,\kappa)$ and $C_2$, such that for all
$k\geq K$,
%
%
\begin{equation}
\label{eqerror-oracle} \sup_{\class_n} L\bigl(\widehat{
\subsp}_m^{\orac}, \widehat{\subsp}_m^{(k),\orac}
\bigr) \leq C_1 M_n\tau_{nm}^2 +
C_2\eps_{nm}^2.
\end{equation}
\end{proposition}
A proof is given in the supplementary material~\cite{supp}, and this completes
step 2.


\subsection{Proof of main results}
\label{secora-act-id}

We now prove the properties of the actual estimating sequence.
The proof relies on the following lemma, which shows the actual and the
oracle sequences are identical up to $2K$ iterations.

%
%
\begin{lemma}
\label{lemmaexclude}
For sufficiently large $n$, 
with probability at least $1-C_0\pn^{-2}$,\vspace*{2pt}
for all $k\leq2K$,
we have
$\widehat{Q}{}^{(k),\orac}_{L\cdot} = 0$,
$\widehat{Q}{}^{(k)} = \widehat{Q}{}^{(k),\orac}$, and hence
$\widehat{\subsp}_m^{(k)} = \widehat{\subsp}_m^{(k),\orac}$.
\end{lemma}
A proof is given in the supplementary material~\cite{supp}, and this completes
step 3.

We now prove Theorems~\ref{thmerror} and~\ref{thmexclude} by
showing that the actual sequence inherits the desired properties from
the oracle sequence. Since Theorem~\ref{thmspecial} is a special case
of Theorem~\ref{thmerror}, we do not give a separate proof.

\begin{pf*}{Proof of Theorem~\ref{thmerror}}
Note that the event on which the conclusions of Lemmas \ref
{lemmaoracle-Sigma}--\ref{lemmaexclude} hold has probability at
least $1-C_0 \pn^{-2}$.
On this event, 
\begin{eqnarray*}
L\bigl(\subsp_m, \widehat{\subsp}_m^{{(\Ks)}}\bigr)
& = & L\bigl(\subsp_m, \widehat{\subsp}_m^{{(\Ks)},\orac}
\bigr)
\\
& \leq& \bigl[ L^{\hf}\bigl(\subsp_m, \subsp_m^{\orac}
\bigr) + L^{\hf}\bigl(\subsp_m^{\orac}, \widehat{
\subsp}_m^{\orac}\bigr) + L^{\hf}\bigl(\widehat{
\subsp}_m^{\orac}, \widehat{\subsp}_m^{{ (\Ks
)},\orac}
\bigr) \bigr]^2
\\
& \leq& C \bigl[ L\bigl(\subsp_m, \subsp_m^{\orac}
\bigr) + L\bigl(\subsp_m^{\orac}, \widehat{
\subsp}_m^{\orac}\bigr) + L\bigl(\widehat{\subsp}_m^{\orac},
\widehat{\subsp}_m^{{ (\Ks
)},\orac}\bigr) \bigr]
\\
& \leq& C_1 M_n\tau_{nm}^2 +
C_2 \eps_{nm}^2.
\end{eqnarray*}
Here, the first equality comes from Lemma~\ref{lemmaexclude}. The
first two inequalities result from the triangle inequality and Jensen's
inequality, respectively. Finally, the last inequality is obtained by
noting that $\Ks\in[K, 2K]$ and by replacing all the error terms by
their corresponding bounds in Lemmas~\ref{lemmaoracle-Sigma}, \ref
{lemmaoracle-S} and Proposition~\ref{proporacle-error}.
\end{pf*}

\begin{pf*}{Proof of Theorem~\ref{thmexclude}}
Again, we consider the event on which the conclusions of Lemmas \ref
{lemmaoracle-Sigma}--\ref{lemmaexclude} hold. Then Lemma \ref
{lemmaexclude} directly leads to the conclusion that $\widehat
{Q}{}^{(k)}_{L\cdot} = \widehat{Q}{}^{(k),\orac}_{L\cdot} =
0$, for all $0\leq k\leq\Ks\leq2K$.
\end{pf*}

\section*{Acknowledgment}

The author would like to thank Iain Johnstone for many helpful
discussions.

\begin{supplement}
\stitle{Supplement to ``Sparse principal component analysis and
iterative thresholding''}
\slink[doi]{10.1214/13-AOS1097SUPP} 
\sdatatype{.pdf}
\sfilename{aos1097\_supp.pdf}
\sdescription{We give in the supplement proofs
to Corollaries~\ref{corrisk} and~\ref{corevec-rate}, Proposition
\ref{prophat-nspike}
and all the claims in Section~\ref{secproof}.}
\end{supplement}

%

\printaddresses

\end{document}